\documentclass[letter,10pt,reqno]{article}
\usepackage{ tikz, pgfplots, todonotes,ulem}
\usetikzlibrary{decorations.markings,arrows, calc}

\usepackage{caption}
\usepackage{amssymb}

\definecolor{shadecolor}{rgb}{0.9, 0.9, 0.81}
\def\Mcal {{\mathcal M}}
\def \eqref #1{ (\ref{#1})}
\def\xt{\tilde{x}}
\def\wh{\widehat}

\def \ds {\displaystyle}

\def\QED{ {\hfill $\blacksquare$}\par \vskip 4pt}
\usepackage{framed,color}

\def\1{\mathbf 1}
\def\wt{\widetilde}

\def \Re { \,\mathrm{Re}\,} 
\def \Im { \,\mathrm{Im }\,} 
\def \z{\zeta}
\def \wh {\widehat}

\newtheorem{definition}{Definition}
\newtheorem{theorem}{Theorem}
\newtheorem{proposition}{Proposition}

\newtheorem{lemma}{Lemma}

\def \bea{\begin{eqnarray}}
\def \pa{\partial}
\def \eea{\end{eqnarray}}
\def \le{\left}
\def \ov{\overline}
\def\ri{\right}

\def\be{\begin{equation}}
\def\ee{\end{equation}}
\def\ben{\begin{displaymath}}
\def\een{\end{displaymath}}
\def\baa{\begin{eqnarray}}
\def\eaa{\end{eqnarray}}

\def\ba{\begin{array}}
\def\ea{\end{array}}
\def\la{\label}
\def\p{\partial}
\def\phi{\varphi}

\makeatletter
\@addtoreset{equation}{section}
\makeatother

\def\Dcal{{\mathcal D}}
\def\pb{{\bf p}}
\def\Rcal{{\mathcal R}}

\def\Lcal{{\mathcal L}}
\def\C{{\mathbb C}}
\def\R{{\mathbb R}}
\def\Z{{\mathbb Z}}
\def\Q{{\mathbb Q}}

\def\Scal{{\mathcal S}}

\def\e{\epsilon}

\def\P{{\mathbb P}}

\def\wh{\widehat}

\def\Pcal{{\mathcal P}}
\def\e{\epsilon}

\def\CP1{{\mathbb C} P^1}
\def\Rcal{{\mathcal R}}
\def\Ccal{\mathcal C}

\def\Qcal{{\mathcal Q}}
\def\Dcal{{\mathcal D}}

\def\db{{\bf d}}

\def\2x2{{\left(\!\!\begin{array}{cc}a&b\\c&d\\\end{array}\!\!\right)}}

\def\f{\frac}
\def\e{\epsilon} 

\def\d{\delta}
\def\d{{ \mathrm d}}
\def\a{\alpha}
\def\b{\beta}

\def\Ch{\widehat{\mathcal C}}

\def\gh{\hat{g}}

\def\dim{{\rm dim}}

\def\gh{\hat{g}}

\def\H1c{H^{(1)}}
\def\H1h{H_{(1)}}
\def\H0c{H^{(0)}}
\def\H0h{H_{(0)}}

\def\a{\alpha}
\def\b{\beta}

\def\gh{{\hat{g}}}

\def\dim{{\rm dim}}
\def\Ocal{{\mathcal O}}

\begin{document}

\title{Discriminant circle bundles over local models of Strebel graphs and Boutroux curves }

\author{M.Bertola and D.Korotkin}

\maketitle

\begin{abstract}
We study 
special "discriminant" circle bundles  over two elementary moduli spaces of meromorphic quadratic differentials with real periods denoted by $\Qcal_0^\R(-7)$ and
 $\Qcal^\R_0([-3]^2)$. The space $\Qcal_0^\R(-7)$      is the moduli space of meromorphic quadratic 
differentials on the Riemann sphere with one pole of order 7 with real periods; it appears naturally   in the study of a neighbourhood of the Witten's cycle $W_1$ in the  combinatorial model  based on  Jenkins-Strebel quadratic differentials of
$\Mcal_{g,n}$. The space  $\Qcal^\R_0([-3]^2)$ 
is the moduli space of meromorphic quadratic 
differentials on the Riemann sphere with two poles of order at most 3 with real periods; it appears in description of a neighbourhood of 
Kontsevich's boundary  $W_{-1,-1}$ of the combinatorial model.
The application of the formalism of the Bergman tau-function to the combinatorial model  (with the goal of computing analytically Poincare dual cycles to certain combinations of tautological classes) requires the study of special sections of circle 
bundles over $\Qcal_0^\R(-7)$ and $\Qcal^\R_0([-3]^2)$; in the case of 
 the space 
$\Qcal_0^\R(-7)$ a section of this circle bundle is given by the argument of the modular discriminant. 
We 
study the spaces $\Qcal_0^\R(-7)$ and $\Qcal^\R_0([-3]^2)$, also called the spaces of Boutroux curves in detail,
together with corresponding circle bundles.

\end{abstract}

\tableofcontents

\section {Motivation and summary of results}

The combinatorial model $\Mcal_{g,n}^{comb}$  of the moduli spaces $\Mcal_{g,n}$ (for $n\geq 1$) based on Jenkins-Strebel differentials \cite{Jenkins,Strebel,Harer,Harer1} was proven to be very fruitful in  studies of geometry of the moduli space.
In particular, it was used in Kontsevich's proof \cite{Kontsevich1, Kontsevich2} of  Witten's conjecture \cite{Witten} about intersection numbers of $\psi$-classes.
 
In this approach the combinatorial model is a CW complex in which each cell is homeomorphic to a simplex labelled by  fatgraphs  embedded
on a Riemann surface.
It was shown in \cite{Igusa2004-1,Igusa2004-2,Mondello1,Mondello2004}, following earlier work \cite{Penner} (where the combinatorial model of $\Mcal_{g,n}$ based on hyperbolic geometry was used)  and \cite{ArbCor} that natural cycles in the combinatorial model
 can be interpreted as Poincar\'e\ duals to certain combinations of Miller-Morita-Mumford classes, or kappa-classes.  These natural cycles consist of cells 
 corresponding to fatgraphs with fixed odd valencies of vertices: cells of highest dimension correspond to three-valent fatgraphs. If the fatgraph contains only one vertex of valency different 
from $3$ then such cycles are called "Witten's cycles"; if the number of 
 such non-generic vertices is bigger than one, the corresponding union of cells is called "Kontsevich cycle". 
  A particular role is played by the following two cycles;
  the first cycle is the Witten's cycle $W_1$ which contains fatgraphs with one vertex of valency $5$ and other vertices of valency $3$.
  The second cycle is the  Kontsevich's boundary $W_{-1,-1}$ of the combinatorial model,  containing fatgraphs which have two vertices of valency $1$  and all other vertices of valency $3$. The Kontsevich's boundary $W_{-1,-1}$ is a subset of
  Deligne-Mumford boundary of $\Mcal_{g,n}$ (some components of the DM boundary are mapped to a point in $W_{-1,-1}$).
  
In general, the Kontsevich--Witten cycle 
$W_{(k_1,\dots, k_r)}$  consists of cells corresponding to Strebel graphs where $r$ vertices have  valency $2k_j+3$ respectively  ($k_j\geq -1$) and all the others have valency $3$ (see \cite{Igusa2004-1}).

 The starting point of computations relating combinatorial cycles to
$\kappa$-classes and their combinations is the explicit expression for the connection forms in  circle bundles (or $U(1)$ bundles) which correspond to tautological line bundles $\Lcal_k$ associated to the  $k$--th puncture, see \cite{Kontsevich1,Zvonkine} and \cite{Igusa2004-1,Mondello1,Mondello2004} for details.

The main motivation of this work is to
provide the  analytical background for an application of the techniques of the  Bergman tau-function  to the combinatorial model $\Mcal^{comb}_{g,n}$ of $\Mcal_{g,n}$. The Bergman tau-function   first appeared in the theory of isomonodromic deformations \cite{Annalen} and
Frobenius manifolds \cite{Dubrovin} in the context on Hurwitz spaces.  It was then extended to the moduli spaces of Abelian \cite{JDG} and quadratic differentials
\cite{Leipzig} over Riemann surfaces. Geometrically, the Bergman tau-function is  a  (holomorphic or meromorphic) section of the product of Hodge line bundle and other  natural line bundles   over the corresponding 
moduli space. Therefore, in a holomorphic framework, an analysis of its singularity structure allows to derive  various relations
in the Picard groups of these moduli spaces
\cite{Advances,MRL,contemp}. 
 The direct application of this technique  in the context of  combinatorial model of $\Mcal_{g,n}$
based on JS differentials is problematic due to the absence of a natural holomorphic structure of the cells.
The   Bergman tau-function appropriately defined  on the Jenkins-Strebel combinatorial model of $\Mcal_{g,n}$ is only
real-analytic in each cell; however, its argument can still be used to get  a smooth section of a circle bundle.
In turn, study of monodromy of these sections  allows us to find cycles of the combinatorial model which are Poincar\'e\ dual to
certain combinations of standard tautological classes on $\Mcal_{g,n}$ \cite{BKfuture}.
The implementation of this idea requires
the study of the local behaviour of the Bergman tau-function and corresponding circle bundle in a neighbourhood of
the cycles $W_1$ and $W_{-1,-1}$.

Hereafter we describe  neighbourhoods of  $W_1 $ and $W_{-1,-1}$ in $\Mcal_{g,n}^{comb}$ 
by applying an appropriate "flat welding" which equivalently can be interpreted via the  "plumbing" construction
using the flat structure introduced on a Riemann surface by the JS differential.

\paragraph{Combinatorial model near  $W_1$ and $W_{-1,-1} $.}
The study of a neighbourhood of the  Witten's cycle $W_1$ in $\Mcal^{comb}_{g,n}$   leads to  the appearance of the real slice 
$\Qcal^\R_0(-7)$ (space of "Boutroux curves") of the moduli space, $\Qcal_0(-7)$,  of meromorphic quadratic differentials 
on Riemann sphere with one pole of order $7$.  

The space $\Qcal^\R_0(-7)$ appears in the study of a neighbourhood of $W_1 \subset\Mcal^{comb}_{g,n} $ as follows;
cells forming  the cycle $W_1 $ are obtained from cells of highest dimension of $\Mcal^{comb}_{g,n}$ by contraction of two edges having exactly 
one  common vertex; this contraction gives rise to the  creation of a ``plumbing zone'' as shown in Fig. \ref{splitW5} which  separates  a component $\C \P^1$ containing the  three  zeros of  $Q$ and a pole of $Q$ 
 of degree $7$  at the nodal point. Thus the quadratic differential $Q$ arising on the separated Riemann sphere belongs to the space $Q_0^\R(-7)$.

The space $\Qcal^\R_0([-3]^2)$ appears
in a neighbourhood of Kontsevich's  boundary $W_{-1,-1} $: it is the moduli space of quadratic differentials on the Riemann sphere with two poles  of order $3$.
The cells forming $W_{-1,-1} $  are obtained from the cells of the highest dimension by simultaneous contraction of two edges having two common 
vertices. Such a contraction creates two ``plumbing zones'' as shown in Fig. \ref{W11-blowup-Fig}; the Riemann sphere arising between the plumbing zones 
contains two simple zeros of $Q$ and two poles of degree $3$ at the arising nodal points. 
Thus the quadratic differential $Q$ arising on the separating Riemann sphere is an element of the space 
$Q_0^\R([-3]^2,[1]^2)$.



 Both spaces $\Qcal^\R_0(-7)$ and $\Qcal^\R_0([-3]^2)$ have real dimension $2$.
Surprisingly enough, we did not find a complete self-contained description of these two elementary moduli spaces in the literature; one of the goals of the paper is to fill this gap.

The stratification of these spaces as well as their complexified versions $\Qcal_0(-7)$ and $\Qcal_0([-3]^2)$, respectively,  is studied in detail in Sections \ref{Q-7}, \ref{Q33}.

\paragraph{Space  $\Qcal_0^\R(-7)$.} 
The generic element of the complex space $\Qcal_0(-7)$  can be represented by the quadratic differential 
\be
Q=(x-x_1)(x-x_2)(x-x_3)(\d x)^2\ ,\ \ \ x_1 + x_2 + x_3=0 \;.
\nonumber 
\ee
The local {\it period} (or {\it homological})  coordinates can be defined as integrals of $v=\sqrt{Q}$ over two independent homology cycles on
the elliptic curve (the "canonical covering")  $\Ch$ defined by the equation $v^2=Q$ (see Sect. \ref{Q-7} for details). The real slice $\Qcal_0^\R(-7)$ of $\Qcal_0(-7)$ is defined by the requirement that all periods of $v$ are real and then the space  $\Qcal_0^\R(-7)$ stratifies  into cells of (real) dimensions $2$ (the top cell of generic elements), of dimension $1$ i.e. the cell where two zeros of $Q$  coincide) and zero dimensional (three coincident roots $x_1=x_2=x_3=0$).

The reality of all periods of $v$ implies that the  period $\sigma$ of the curve $\Ch$ belongs to a one-dimensional subset $\Rcal_1$ of the modular curve shown in Fig.\ref{sigma71int}, left pane. 

The space  $\Qcal^\R_0(-7,[1]^3)$ is fibered over the set $\Rcal_1$ 
with  fiber  $\R_+$.   The points at $\infty$ of $\Rcal_1$ correspond to  the space  $\Qcal^\R_0(-7,1,2)$. The point of intersection of $\Rcal _1$ with the real  axis in the plane of the $J$--invariant (Fig.\ref{sigma71int}, right pane) is $J\simeq 940.34 $ and it corresponds to the Boutroux--Krichever curve  \cite{Wiegman, Krichever, Bertola}, which is the unique Boutroux curve in   $\Qcal^\R_0(-7,[1]^3)$ possessing a real involution $^\star$ which leaves the Jenkins-Strebel differential invariant: $\ov{Q(x^\star)} = Q(x)$. 

The zeros $x_1$, $x_2$ and $x_3$ are connected by two horizontal (i.e. parallel to the real line in the plane of coordinate $z(x)=\int^x v$) geodesics in the metric $|Q|$. We will always denote by $x_2$ the "central" zero which is connected by these geodesics to two others. Then the remaining zeros $x_1$ and $x_3$ 
can also be labeled such that in the positive direction around the origin $x_3$ goes after $x_2$ and $x_1$ after $x_3$.

The lengths of the two geodesics are given by the absolute values  $A = \Big|\int_{x_1}^{x_2} \sqrt{Q}\Big|, \ B=\Big|\int_{x_2}^{x_3} \sqrt{Q}\Big|$.  The lengths $(A,B)$  define the map of the space $\Qcal^\R_0(-7,[1]^3)$ to $\R_+^2$.

\paragraph{Discriminant circle bundle on  $\Qcal_0^\R(-7)$.}  Consider the section $\phi_1$ of a $U(1)$ bundle over $\Qcal^\R_0(-7)$ given by the argument of the modular discriminant
\be
\phi_1={\rm Arg} \Delta_1 
\ee
where 
\be
\label{Delta5}
\Delta_1 =[(x_1-x_2)(x_3-x_1)(x_3-x_2)]^2.
\ee
Although $\Delta_1$ itself is well-defined only on the stratum  
$\Qcal^\R_0([1]^3,-7)$, its argument $\phi_1$ can  be  continued smoothly through the boundaries $A=0$, $B=0$.

We will show (Theorem \ref{propdelta5})  that then increment   of $\phi_1$ on $\Qcal^\R_0(-7)$  between the boundaries $A=0$ and  $B=0$ equals $\pi/5$; this computation is technically non-trivial. Therefore, $\pi/5$ is the monodromy of $\phi_1$ around the
 zero dimensional cell (represented by $Q = x^3\d x^2$) on the space $\Qcal^\R_0(3,-7)$.

\vskip0.5cm
 
\paragraph{Description of $\Qcal_0^\R([-3]^2)$.} 
The complex moduli space 
$\Qcal_0([-3]^2)$ consists of   quadratic differentials $Q$ on the Riemann sphere with  two poles of degree $3$. 
 The generic element of this space is  represented by 
\be
Q=\frac{(x-x_1)(x-x_2)}{x^3}(dx)^2\ ,\ \ \ x_1\neq x_2,\ \ x_1\neq 0 \neq x_2.
\ee
The  integrals of $\sqrt{Q}$ along paths two arbitrary independent homology classes on the
canonical covering $y^2=x(x-x_1)(x-x_2)$ are  local period coordinates  on  $\Qcal_0([-3]^2)$.

The space  $\Qcal^\R_0([-3]^2)$ is defined by the condition that all periods of $v=\sqrt{Q}$ on $\Ch$ are real; it is fibered over the  set  $\Rcal_{-1,-1}$ of real dimension $1$ within the modular curve (Fig. \ref{setR11}, left pane)  with the fiber $\R_+$.

The points at $\infty$ of $\Rcal_{-1,-1}$ correspond to the ``wall'' i.e. to the space  $\Qcal^\R_0(2,[-3]^2)$. 

The set  $\Rcal _{-1,-1}$ intersects the real  axis in the plane of the $J$--invariant at $J\simeq \{-1690, 586, 7791\}$ (Fig.\ref{setR11}, right pane). 
The points $J\simeq -1690$ and $J \simeq 586$ correspond to curves with a real involution, however the differential $Q$ is not invariant with respect to this involution.  
The point  $J\simeq 7791 \ (x_1\simeq  1.8037, \ \ x_2 = -0.3797)$ corresponds to a  curve  possessing a real involution $\star$ (acting as a standard complex conjugation) such that $\ov{Q(x^\star)} = Q(x)$. This is the natural analog of the Boutroux--Krichever curve in the space $\Qcal^\R_0([-3]^2)$. 

One of the zeros, which we denote by $x_1$, is connected by an infinite horizontal geodesics in the metric $|Q|$ to the pole $x=0$, and another zeros ($x_2$) to the pole $x=\infty$. In this way we get the labeling of the zeros $x_1$ and $x_2$ in this case.
The zeros $x_1$ and $x_2$ are connected by two horizontal geodesics which we denote by $e_1$ and $e_2$.
The geodesics $e_1$ and $e_2$ will be labelled as follows: when one arrives to $x_1$ along horizontal geodesics emanating from $x=0$, turning
right one follows the geodesics $e_1$, and turning left one follows the geodesics $e_2$ (see Fig. \ref{fig2degmid}). The length  of $e_1$ will be denoted by $A$ and the length of $e_2$ by $B$.

The lengths $(A,B)$ 
define the map of the space $\Qcal^\R_0([-3]^2,[1]^2)$ to $\R_+^2$.

\paragraph{Discriminant circle bundle on $\Qcal_0^{\R}([-3]^2)$.} The  sections $\phi_{-1,-1}^\pm$ of the two natural $U(1)$ bundles over $ \Qcal_0^\R([-3]^2) $  are given by the following expressions:
\bea
&&\phi_{-1,-1}^\pm = {\rm Arg}\,  \Delta^\pm_{-1,-1}
\eea
\bea
&& \Delta^+ _{-1,-1} := x_1^6 x_2^6 (x_1-x_2)^2,\  \ \qquad  \Delta^- _{-1,-1} := x_1^6 x_2^6 (x_1-x_2)^{26}\;.
\eea
The  analysis contained in Section \ref{space11}  shows that the increments of $\phi^{\pm}_{-1,-1}$ 
between the boundary $A=0$ and the boundary $B=0$ on $Q^\R_0([-3]^2)$ equals $ {13\pi} $ and $  { 25\pi} $, respectively.

Therefore, monodromies of $\phi^{\pm}_{-1,-1}$ 
along a simple closed, non-contractible  loop in ${\Qcal}^\R_0([-3]^2)$ 
equal ${13\pi} $ and $  { 25\pi} $, respectively  (see Thm. \ref{prop2}).

\vskip 4pt
In summary, the main results of this paper are the following.  First,  we show how the moduli spaces  $Q^\R_0(-7)$ and $Q^\R_0([-3]^2)$
  appear under a degeneration of a generic Strebel graph via shrinking of two adjacent edges. Second, we study 
analytically the spaces   $Q^\R_0(-7)$ and $Q^\R_0([-3]^2)$ and natural circle bundles over them. These results provide the analytical tools necessary to
apply the formalism of 
the Bergman tau-function in the description of various tautological classes to the combinatorial model of  $\Mcal_{g,n}$ \cite{BKfuture}.

The paper is organized as follows.
In Section \ref{combmod} we recall the combinatorial model of $\Mcal_{g,n}$ based on Jenkins-Strebel differentials.
In section \ref{localW5} we describe a neighbourhood of Witten's cycle $W_1 $ in the combinatorial model and show how the 
moduli space $Q_0^\R(-7)$ appears in this context. In Section \ref{localW11} we describe a neighbourhood of Kontsevich's boundary $W_{-1,-1} $ of $\Mcal_{g,n}^{comb}$ and demonstrate   the appearance of the space $Q_0^\R([-3]^2)$.
In Section \ref{space5} we study the geometry of the space $Q_0^\R(-7)$ in detail and compute the increment of the argument of
the modular discriminant $\Delta_1 $ on this space.
Finally, in Section \ref{space11} we study the geometry of the space $Q_0^\R([-3]^2)$ and compute the increment of the
arguments of $\Delta^\pm_{-1,-1} $ on this space.

\section{Combinatorial model of $\Mcal_{g,n}$ via Jenkins-Strebel differentials}
\la{combmod}

Here we briefly recall the main ingredients of the combinatorial model of the moduli spaces of Riemann surfaces
based on Jenkins-Strebel differentials. 
The  moduli space of Riemann surfaces of genus $g$ with $n$ marked points is denoted by $\Mcal_{g,n}$ and its Deligne-Mumford compactification by $\overline{\Mcal}_{g,n}$. Let $\Ccal$ be a Riemann surface of genus $g$ and $Q$ be a meromorphic quadratic differential with second order poles at  the
points $y_1,\dots,y_n$. 
Zeros of  $Q$ are denoted by $x_1,\dots, x_m$ and their multiplicities by ${\bf d}=(d_1,\dots,d_m)$;
we have $\sum_{i=1}^m d_i= 4g-4+2n$.
 Denote by $\Qcal_{g,n}^{{\bf d}}$ the moduli space of such quadratic differentials; its dimension equals
$2g-2+n+m$. For the stratum of highest dimension, when all zeros are simple and $m=4g-4+2n$,  this dimension  equals $6g-6+3n$. 

\paragraph{Canonical cover.} The equation 
$$v^2=Q$$
in $T^*\Ccal$ defines a two-sheeted covering $\Ch$ of $\Ccal$, called {\it canonical covering} in Teichm\"uller theory
($\Ch$ is known under the name of "spectral covering" in the theory of Hitchin's systems or under  the name of "Seiberg-Witten curve" in the theory of supersymmetric Yang-Mills theory). The branch points of the covering coincide with zeros 
of odd multiplicity  of $Q$, whose  number we denote by $m_{odd}$:  then the genus of the canonical covering is $\gh=2g+\f{m_{odd}}{2 }-1$.  For a  generic $Q$, which has
$4g-4+2n$ simple zeros, $\gh=4g-3+ n$. It is convenient to introduce the notation
$$
g_-=\gh-g= g+\frac{m_{odd}}{2}-1
$$

Denote by $\mu$ the holomorphic involution interchanging the sheets of $\Ch$. 
The Abelian differential of the third kind $v$ on $\Ch$ has   $2n$ simple poles at the points $\{y_i,y_i^\mu\}$
(slightly abusing the notations we use $y_i$ to denote positions of poles both on $\Ccal$ and $\Ch$).
Consider the decomposition of $H_1(\Ch\setminus \{y_i,y_i^\mu\}_{i=1}^n ,\R)$ (whose dimension equals $2\gh+2n-1$) into the direct sum of even and odd subspaces:
\be
H_+\oplus H_-\ ,\qquad \dim H_+=2g+n-1,\ \ \ \dim H_-=2 g_-+n.
\ee
 Notice that
$
\dim H_-=\dim \Qcal_{g,n}^{{\bf d}}\;.
$
Consequently, one can introduce a  system of local coordinates on the  space  $\Qcal_{g,n}^{{\bf d}}$ called {\it ``homological coordinates''}  by choosing a set  of independent cycles $s_i$, $i=1,\dots,\dim H_-$ in $H_-$ and
integrating $v$ over these cycles:
\be
\Pcal_i=\int_{s_i} v\;.
\ee
The following definition of Jenkins-Strebel differential is  equivalent to the standard one \cite{Strebel}:
\begin{definition}
The quadratic differential $Q$ is called Jenkins-Strebel differential if all  homological coordinates $\Pcal_i$ are real.
\end{definition}

The reality of all $\Pcal_i$ implies that the biresidues of $Q$ at its poles $y_j$ are real and negative i.e. there exist such $p_j\in \R_+$ 
that in any local coordinate $\zeta$ near $y_j$ one has:
\be
Q(\zeta) =\frac {-(p_j/2\pi) ^2}{\zeta^2}\le(1 + \mathcal O(\zeta)\ri) (\d \zeta)^2
\la{locQ}
\ee

Results of Jenking and Strebel provide the   existence and uniqueness of a differential $Q$ on a given Riemann surface with $n$ marked points with {\it given}
constants $p_i\in\R_+$ and all {\it real} periods $\int_{s_i}v$ for $s_i\in H_-$ (combinations of cycles surrounding  $z_i$ also form a part of $H_-$; the  reality of periods of $v$ 
around these cycles is guaranteed by the form of biresidues in (\ref{locQ})). This statement provides a basis for the combinatorial description of $\Mcal_{g,n}$.

For each given vector  $\pb=(p_1,\dots, p_n) \in \R_+^n$ one can construct a  combinatorial model of the moduli space $\Mcal_{g,n}$ which
will be denoted by  $\Mcal_{g,n}[\pb]$.
Namely, for each Riemann surface $\Ccal$ of genus $g$ with $n$ marked points $z_1,\dots,z_n$ consider the unique Jenkins-Strebel differential $Q$ 
whose singular part at $z_i$  is as in (\ref{locQ}). 

The stratum  $\Mcal_{g,n}^{\db}[\pb]$ of  $\Mcal_{g,n}[\pb]$ consists of punctured Riemann surfaces  such that the Jenkins-Strebel differential
$Q$ has zeros of multiplicities ${\bf d}=(d_1,\dots,d_m)$. The real dimension of the stratum   $\Mcal_{g,n}^{\db}[\pb]$ is given by
\be
\dim_\R \Mcal_{g,n}^{\db}[\pb]= 2g-2+m
\ee
The  largest stratum corresponds to Jenkins-Strebel differentials with simple zeros and has real dimension 
equal to  $6g-6+2n$ (in this case $m=4g-4+2n$)  which coincides with the real dimension of $\Mcal_{g,n}\times \R_+^n$.

The oriented horizontal trajectories of the 1-form $v$ on $\Ch$ (which project down to the non-oriented horizontal trajectories of $Q$ on $\Ccal$) 
connect zeros $x_i$.

These horizontal trajectories on $\Ccal$ form the edges of an embedded graph (also called {\it fatgraph}, or {\it ribbon graph}) $\Gamma$  on $\Ccal$ with vertices at $x_1,\dots,x_m$ and
valencies $d_1+2,\dots,d_m+2$, respectively. Therefore, the stratum of highest dimension $6g-6+2n$ corresponds to fatgraphs on $\Ccal$  with trivalent vertices.

\begin{figure}
\begin{center}
\includegraphics{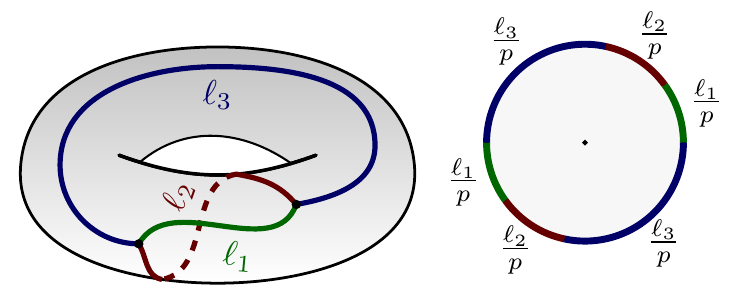}
\end{center}
\caption{Fatgraph on a genus $1$ Riemann surface representing a point  in $\Mcal_{1,1}[\pb]$. Two simple zeros $x_1$ and $x_2$ of $Q$ are connected by 3 edges of lengths $\ell_1$, $\ell_2$ and $\ell_3$. The fatgraph has only one face of perimeter $p=2(\ell_1+\ell_2+\ell_3)$. The lengths $\ell_1$ and $\ell_2$ can be used as coordinates on  $\Mcal_{1,1}[p]$}
\label{M11}
\end{figure}

The length $\ell_e$ of an edge $e\in E(\Gamma)$ connecting two vertices $v_1, v_2$  is  equal to the absolute value of the integral of $v$ along the horizontal trajectory connecting the two vertices. In turn, up to  a factor of $1/2$,  this length coincides with the integral of $v$ over the integer cycles in  $H_-$ consisting of the trajectory $e$ on one sheet of the projection  $\Ch \to \Ccal$ and the same trajectory in the opposite direction on the other sheet.
Fixing the vector $\pb\in \R_+^n$ one imposes $n$ linear constraints on the  lengths of the edges: for  the  $j$--th face $f_j \in F(\Gamma),\ \ j=1,\dots, n$ we have $\sum_{e\in \pa f_j} \ell_e = p_j$. 
A simple example of fatgraph on a  Riemann surface of genus $1$ with one marked point is shown in Fig. \ref{M11}.
In the flat metric $|Q|$ on $\Ccal$ the neighbourhoods of a pole $y_i$ is an infinite cylinder of perimeter  $p_i$.\par\vskip 4pt

The union of all strata  $\Mcal_{g,n}^{\db}[\pb]$ for fixed $\pb\in \R_+^n$ forms the combinatorial model  $\Mcal_{g,n}[\pb]$ of $\Mcal_{g,n}$. This combinatorial 
model is  set-theoretically isomorphic to $\Mcal_{g,n}$.

The compactification  ${\ov {\Mcal}}_{g,n} [\pb]$ of the combinatorial model is constructed by addition of the so--called {\it Kontsevich boundary}  to $\Mcal_{g,n} [\pb]$: 
the stratum $W_{-1,-1}$  of real co-dimension $2$ of the Kontsevich boundary corresponds to fatgraphs with exactly two one--valent vertices (i.e. two simple poles of the JS differential $Q$), while all other vertices remain trivalent.

The Kontsevich boundary is  ``smaller'' than the Deligne-Mumford boundary of $\Mcal_{g,n}$ since  the Jenkins-Strebel  combinatorial model  
is not applicable to non-punctured Riemann surfaces.

Each face $F_j$ of the fatgraph $\Gamma$ (the face $F_j$ contains the pole $y_j$) can be mapped to the unit disk via the map
\be
w_j(x)=\exp\left(\frac{2\pi i}{p_j} \int_{x_j}^x v\right)\;;\
\la{faceunit}
\ee
where $x_j$ is an arbitrarily chosen vertex of the face $F_j$. Under the map (\ref{faceunit}) $y_j$ is being mapped to the origin and $x_j$ to 1.
One can always choose the branch of the differential $v=\sqrt{Q}$ which has  residue $\frac{p_j}{2\pi i} $ at $y_j$. 
Within  the face $F_j$ the flat  metric $|Q|=|v|^2$ on $\Ccal$    is expressed as:
\be
|Q|= \frac{p_j^2}{4\pi^2} \left|\frac{dw_j}{w_j}\right|^2\;.
\ee

The fatgraph corresponding to the stratum of highest dimension in the combinatorial model $\Mcal_{1,1}^{comb}[p]$ shown in Fig. \ref{M11} 
has only one face; this face can be mapped to the unit disk. The constraint between lengths  in this case reads $2(\ell_1+\ell_2+\ell_3)=p$.

Inverting this logic  it is possible to construct a polyhedral Riemann surface (i.e. Riemann surface with flat metric with conical singularities)
from a fatgraph equipped with lengths of all edges by the {\it conformal welding} \cite{Strebel}.

All strata of $\Mcal_{g,n}[\pb]$ can be obtained from the stratum of highest dimension $W$ (which corresponds to fatgraphs with all tri-valent  vertices)
by contraction of one or more edges. Various components of $W$ are glued along (real) codimension one boundaries which correspond to 
 fatgraphs with one vertex of valency $4$ and all other vertices of valency $3$ (i.e. the 
Jenkins-Strebel differential corresponding to the boundary  has one double zero and other zeros are simple).
The procedure of crossing the boundary  from one cell of $W$ to another 
(sometimes in this way one can get into the same cell; then the cell is "wrapped on itself") is described by the so-called Whitehead move shown in Fig.\ref{flip}.

 Our goal will be to study in detail the contraction of  two edges having one vertex in common. Then depending on geometry  one gets  either  cells forming
the Witten's cycle $W_1 $ (corresponding fatgraphs have one 5-valent vertex and other vertices of valency $3$), or cells forming the  Kontsevich's boundary 
$W_{-1,-1} $ of $\Mcal_{g,n}[\pb]$ (which corresponds to fatgraphs having two vertices of valency 1 and all other vertices of valency $3$).
These two types of contraction are shown in Fig.\ref{figw5deg}

\noindent
\begin{minipage}{0.6\textwidth}
\begin{center}
\includegraphics[width=1\textwidth]{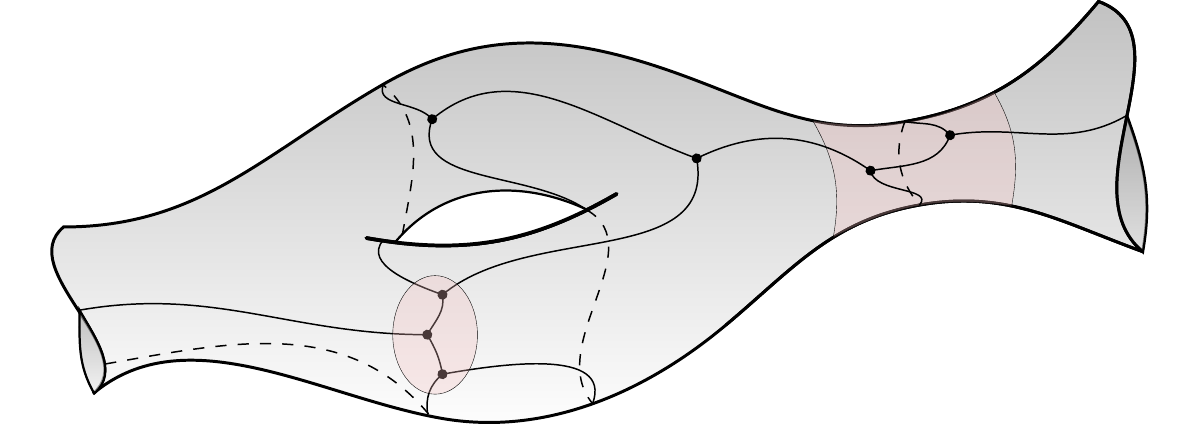}
\end{center}
\captionof{figure} {A typical embedded Strebel graph. The two highlighted  regions are where the contraction of edges leads either to the
cell $W_1 $ (left region) or to the cell $W_{-1,-1} $ (right region). }
\label{figw5deg}
\end{minipage}
\begin{minipage}{0.4\textwidth}
\begin{center}
\includegraphics{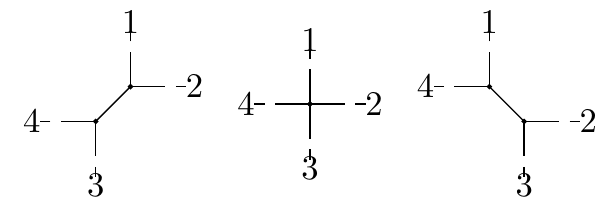}
\end{center}
\captionof{figure}{The Whitehead move.}
\label{flip}
\end{minipage}

\section{ A local model near Witten's cycle $W_1$: the  space $Q_0^\R(-7)$}
\la{localW5}

Cells of the  cycle  $W_1$ in the combinatorial model $\Mcal_{g,n}[\pb]$ 
(the vector $\pb$ will be kept fixed in all constructions below) 
correspond to fatgraphs which have all vertices of valency 3 except one vertex which has valency 5. 
Here we shall understand an $n$-punctured   Riemann surface $\Ccal$ as an element of the combinatorial model 
$\Mcal_{g,n}[\pb]$ i.e. we assume that $\Ccal$ uniquely defines the JS differential $Q$.

A point of $W_1 $ can be obtained from a trivalent fatgraph  via contraction of two 
edges having one vertex in common; conversely, any Riemann surface from a neighbourhood of $\Ccal$ in $W$ is  a two (real)-parameter deformation of $\Ccal$. We are going to show how this  deformation can be described via "flat surgery"  of a Riemann surface $\Ccal$ 
and a Riemann sphere equipped with an appropriate real-normalized quadratic differential.
This flat surgery also explicitly gives the Jenkins-Strebel differential and the fatgraph  on the deformed Riemann surface.

Furthermore, the flat surgery can also be understood in terms of  a  "plumbing construction" connecting $\Ccal$ with a Riemann sphere equipped with an appropriate flat metric.

Denote the zero of order $3$ of the JS differential on $\Ccal$ by $x_1$; denote also the  fatgraph corresponding to $\Ccal$ by $\Gamma$. Let $\zeta$ be a ``distinguished'' local parameter on $\Ccal$ near $x_1$: $\zeta(x)=\left[\frac{5}{2}\int_{x_1}^x v\right]^{2/5}$ ($\zeta(x)$ is defined up to a fifth root of unity); thus $Q(x)$ can be written in terms of $\zeta$ as $Q(\zeta)=\zeta^3(\d\zeta)^2$.  The local coordinate on the canonical cover $\Ch$ near $x_1$ is given by $\zeta^{1/2}$.

Introduce  also the ``flat'' coordinate (both on $\Ccal$ and $\Ch$):
$$
z(x)=\int_{x_1}^x v\;;
$$
in a neighbourhood of $x_1$ the flat coordinate is expressed via the local coordinate $\zeta$  as $z=\frac{2}{5} \zeta^{5/2}$. We have $Q(x) = (\d z(x))^2$ on $\Ccal$ and $v=\d z$ on $\Ch$. The metric $|Q(x)|$ 
has a conical point at $x_1$ with cone angle $5\pi$ while at all other vertices the cone angle equals $3\pi$.

\subsection{Resolution of a $5$-valent vertex by flat surgery}
\label{res5}

\begin{figure}[htb]
\begin{center}
\includegraphics{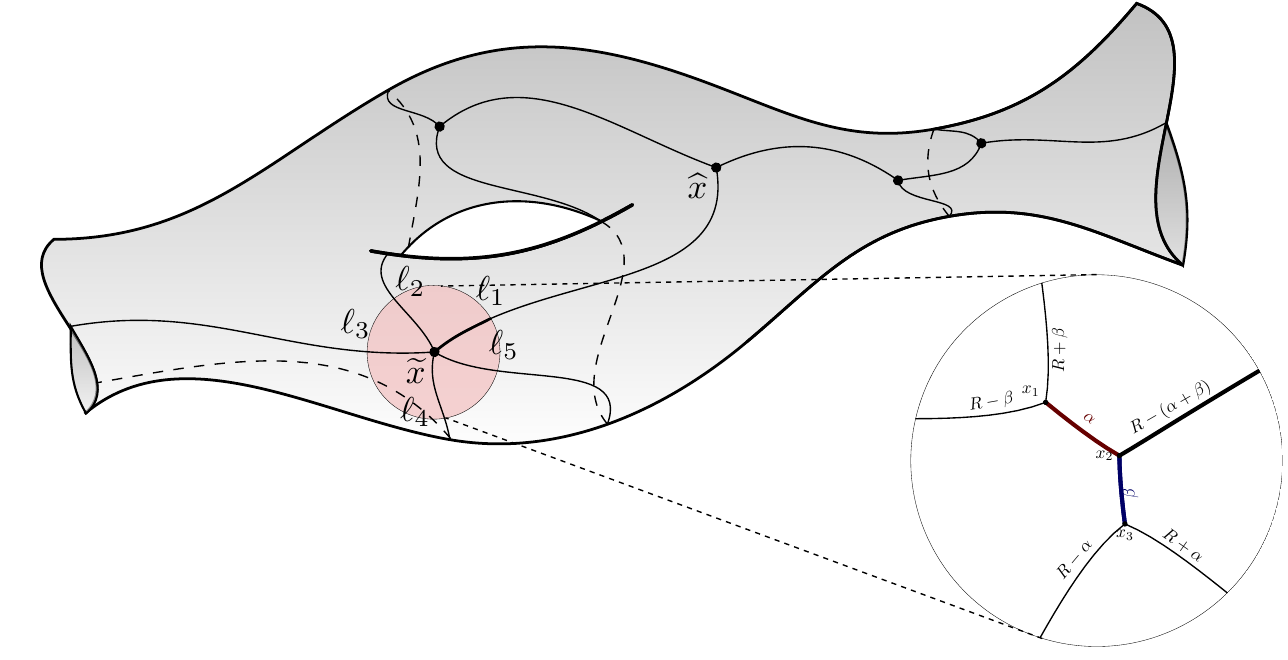}
\end{center}

\captionof{figure}{The Riemann surface $\mathcal C^{\a,\b}$ is obtained by replacing the five-valent vertex with three regular vertices and two new edges of length $\a$ and $\b$ between them with appropriate adjustment of the lengths of the original fatgraph. The limit $\a,\b \to 0 $ produces the original Riemann surface.
The  five arcs on the boundaries of the plug-in region are five semicircles or radius $R$ in the flat coordinates of either metrics.}
\label{W5plugin}
\end{figure}

There are $5$ edges of the fatgraph $\Gamma$ emanating from the vertex $x_1$; their directions are given by the angles $\frac{2\pi k}{5}$ in the $\zeta$-plane. Denote their lengths (starting from the edge going along positive real line) by $\ell_1,\dots,\ell_5$.
We are going to construct an explicit deformation $\Ccal^{\alpha,\beta}$ of $\Ccal$ with two real ("small") parameters, $\alpha$ and $\beta$, by changing these lengths as follows:
\be
(\ell_1,\dots,\ell_5)\to (\ell_1-(\alpha+\beta),\, \ell_2+\beta,\, \ell_3-\beta,\, \ell_4-\alpha, \,\ell_5+\alpha)
\la{deforl}\ee
The triple zero, $\wt x$, of $Q$ on $\Ccal$ splits into three simple zeros $(x_1,x_2,x_3)$ 
of the JS differential $Q^{\a,\b}$ on $\Ccal^{\a,\b}$. We assume that the new edges of length $\ell_1-(\alpha+\beta)$ and $ \ell_2+\beta$ meet at 
$x_1$; $x_1$ is connected to $x_2$ by the edge of length $\a$, $x_2$ is the endpoint of the edge $\ell_2+\beta$ and it is also connected to $x_3$ by the edge of length $\beta$; the edges $\ell_4-\alpha$ and $\ell_5+\alpha$ meet at $x_3$. The resulting
configuration is shown in Fig.\ref{W5plugin}.
\paragraph{}
We then construct a two--(real) parameter family  $\Ccal^{\a,\b}$ of  deformations of $\Ccal$  as explained below.

Fix some radius $R$ which satisfies the conditions $R<\ell_i$ for $i=1,\dots,5$. 
 Consider the disk $\Dcal_R$ given by  $|z(x)|\leq R$ on $\Ccal$ (see Fig. \ref{W5plugin}); in terms of the coordinate $\zeta$ the disk $\Dcal_R$ is given by $|\zeta|\leq \left(\frac{5R}{2}\right)^{2/5}$.
Five edges of the fatgraph $\Gamma$ within $\Dcal_R$ are given by the segments  
$ \le[0,\le(\frac {5R}2\ri)^\frac 25 \ri ]{\rm e}^{\frac{2i\pi}5 k}$, $k=0,\dots, 4$ in the $\zeta$-plane.
The perimeter of $\Dcal_R$ in the flat metric $|Q|$ equals $5\pi R$.

Let us now excise the region $\Dcal_R$ from the Riemann surface $\Ccal$  and obtain a Riemann surface 
with boundary, which we denote by $\Ccal_R$.  Denote by $\Ch_R$ the canonical cover $\Ch$ with the disk $\Dcal_R$ deleted on both copies of $\Ccal$.
Let us assume that the branch cut connecting the  triple zero $\wt x$ with some other zero on $\Ccal$ (say, $\wh x$) goes along the positive real line in the $\zeta$ - coordinate.

We are going to attach to $\Ch_R$ another disk $\Dcal^{\a,\b}_R$ which we now describe.
Consider an element of the moduli space $Q_0^{\R}[-7]$, that is, a quadratic differential $Q_0$ on $\C P^1$ given by
\be
\label{Q0}
Q_0(x) =(x-x_1)(x-x_2)(x-x_3)  \d x^2 
\ee 
such that $x_1+x_2+x_3=0$.
The Abelian differential $v_0=\sqrt{Q_0}$ is defined on the canonical cover $\Ch_0$ of $\CP1$ that is the elliptic 
curve with branch points at $x_1,x_2,x_3$ and $\infty$.

Assume that all periods of $v_0$ on $\Ch_0$ are real. Then there exist two horizontal trajectories of $v_0$ on $\CP1$ which connect $x_1$, $x_2$ and $x_3$. Assume that $x_2$ is the ``central'' zero i.e. it is connected by the horizontal trajectories to $x_1$ and $x_3$. Choose the branch cuts
on $\Ch_0$ to connect $x_2$ with $x_1$, $x_3$ and $\infty$ along the horizontal trajectories. Assume also that $x_1$ and $x_3$ are enumerated such tat
the branch cuts $[x_2,\infty)$, $[x_2,x_1]$ and $[x_2,x_3]$  meet at $x_2$ in the counterclockwise order as shown in Fig.\ref{localmodelW5}, left (this picture is drawn in the coordinate $x$).

Introduce on $\Ch_0$ a canonical basis of cycles $(a,b)$ such that the $a$-cycle encircles the branch points $x_1$ and $x_2$ and $b$-cycle encircles the branch points $x_2$ and $x_3$. The real periods of $v_0$ over cycles $a$ and $b$ are expressed as follows via the lengths $\a$ and $\b$
of the branch cuts $[x_2,x_1]$ and $[x_2,x_3]$:
\be
\Big|\int_a v_0\Big|=2\a\;\;,\hskip0.7cm 
\Big|\int_b v_0\Big|=2\b\;.
\la{intv0}
\ee
The conditions (\ref{intv0}) determine the branch points $x_i$ up to multiplication of all $x_i$ by a fifth root of unity (see Th. \ref{thm123}).

The quadratic differential $Q_0$ which has three simple zeros and one pole of degree 7  on $\CP1$ 
is an element of the space $\Qcal_0^{\R}([1]^3,-7)$; it is an analog of the 
Jenkins-Strebel differential in the case of a pole having an order higher than 2.
\vskip0.5cm
\noindent\begin{minipage}{0.6\textwidth}
\includegraphics[width=0.99\textwidth]{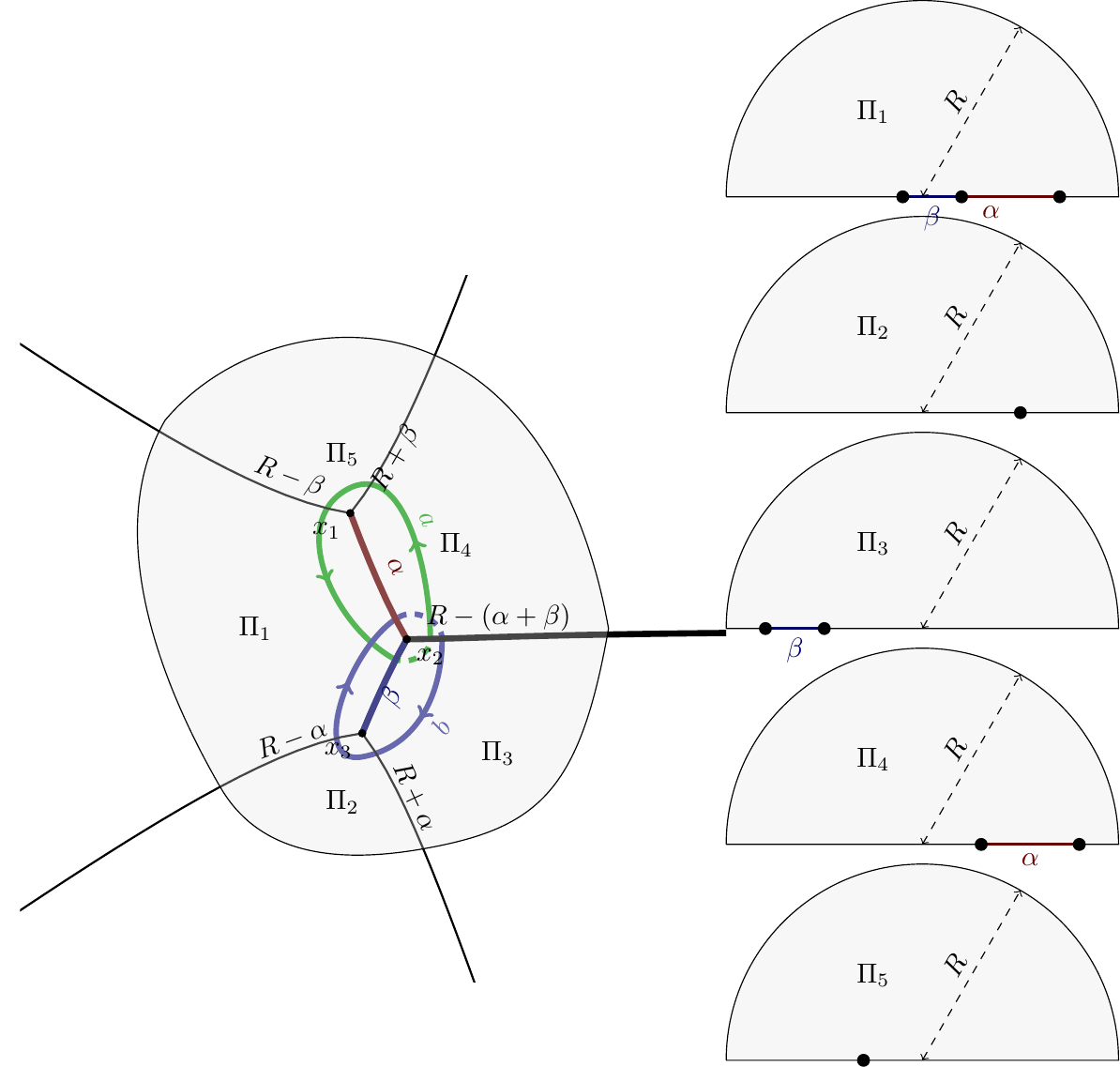}
\end{minipage}
\hspace{0.01\textwidth}
\begin{minipage}{0.39\textwidth}
\captionof{figure}{ Left: the ribbon graph on $\C\P^1$ of the quadratic differential $Q_0$ and the geodesic disk of radius $R$ drawn in coordinate $x$. Indicated are the distances in the metric $|Q_0|$ along the trajectories. Right: the five half-disks in the upper half plane, uniformized by the flat coordinate.}
\label{localmodelW5}
\end{minipage}
\vskip0.5cm
We choose the branch cuts  to go along the edges $[x_1, x_2]$, $[x_3,x_2]$ and $[x_2,\infty]$ (the thick lines in 
Fig. \ref{localmodelW5}, left pane).
The edges of the fatgraph $\Gamma_0$ split the $x$-plane into 5 regions
 as shown in Fig.\ref{localmodelW5} (left). Five horizontal trajectories of $v_0$ which connect $x_1$, $x_2$ and $x_3$ to $\infty$ approach $\infty$ along the rays ${\rm arg}\, x= 2\pi k/5$. 
In the flat metric $|Q_0|$ on $\CP1$ each of these five regions is uniformized by the flat coordinate $z$ to a half-plane $\Im z>0$. 

  On each of the five rays  connecting the point at infinity with the points $x_i$ we cut  segments whose lengths (starting from the ray
$[x_2,\infty)$ and counting couterclockwise)  are chosen to be  $R-\a-\b$, $R+\beta$, $R-\beta$, $R-\a$ and $R+\a$.
Then we consider a half-circle (in the flat metric $|Q_0|$) in each region which connects the endpoints of the chosen segments 
(Fig. \ref{localmodelW5}, right). It is easy to see that the diameters of all of these five half-circles coincide and are all equal to $2R$.
In this way we obtain five regions which are half-disks of radius $R$ in the flat coordinate $z$; they will be denoted by $\Pi_1,\dots,\Pi_5$, and their union by $\Dcal_R^{\a,\b}$.

The flat coordinates in each of the regions $\Pi_i$ are defined up to sign and translations; the  global  flat coordinate on the whole $\Dcal_R^{\a,\b}$  is defined as follows. Choose the initial point of integration to be $x_2$ and choose the system of three branch cuts within $\Dcal_R^{\a,\b}$ as shown in Fig.\ref{localmodelW5} (left):  then the flat coordinate 
on  $\Dcal_R^{\a,\b}$  with three deleted branch cuts is defined by
\be
w(x)=\a+\b+\int_{x_2}^x v_0 
\la{flatcoordab}
\ee
where the determination of $v_0$ has to be chosen so that the image of $\Pi_4$ is in the upper half plane of the $z$--variable.
One can easily verify that,  according to the orientation of the $a$- and $b$- cycles on $\Ccal_0$ shown in Fig.\ref{localmodelW5}, left pane, the values of the coordinate $w(x)$ on different 
sides of the branch cut $[x_2,\infty)$ differ by a sign.

The canonical two-sheeted cover of the domain  $\Dcal_R^{\a,\b}$ with branch cuts going along horizontal trajectories $[x_2,x_1]$, $[x_2,x_3]$ and $[x_2,\infty)$ 
will be denoted by $\widehat{\Dcal}_R^{\a,\b}$.

The shape and perimeter of the boundary of $\Dcal_R^{\a,\b}$ in the global flat coordinate is a circle winding by an angle of $5\pi$. The  positions of the points where the  outgoing edges (in $z$-coordinate) 
and the branch cut
intersect the boundary of $\Dcal_R^{\a,\b}$, coincide with those of $\Dcal_R$. In particular, the perimeter of $\Dcal_R^{\a,\b}$ also equals $5\pi R$.
  Therefore, we can identify the boundary
of $\Ccal_R$ with the boundary of $\Dcal_R^{\a,\b}$ (this procedure is called "conformal welding" in \cite{Strebel})
to get the new Riemann surface $\Ccal^{\a,\b}$ together with its canonical cover $\Ch^{\a,\b}$.  The five lengths of the 
edges originally connected to the zero of order 3 on $\Ccal$ are then changed according to (\ref{deforl}). Notice that this deformation preserves perimeters of all faces of the fatgraph; therefore the Riemann surface $\Ccal^{\a,\b}$ and the corresponding fatgraph 
belong to the same combinatorial model $\Mcal_{g,n}[\pb]$ as the original surface $\Ccal$.

Moreover, the result of such "conformal welding" does not depend on the choice of radius $R$ of the excised disk as long as it remains sufficiently small in comparison with $\ell_1,\dots, \ell_5$ and sufficiently large in comparison with $\a$ and $\b$.

 \subsection{Plumbing construction}

To study the limit $\a,\b\to 0$ in the previous "conformal welding" scheme one needs to degenerate the quadratic differential 
(\ref{Q0}) on $\CP1$. The equivalent "plumbing" construction presented here allows to keep $\Ccal$ and $\Ccal_0$ fixed.

To  deform $\Ccal$ into  $\Ccal^{\a,\b}$ (and conversely, to study the limit $\a,\b\to 0$ of $\Ccal^{\a,\b}$)
we introduce the "plumbing parameter" $t=\a+\b$ and the parameters
\be
A=\f{\a}{t}=\f{\a}{\a+\b}\;\;\hskip0.7cm
B=\f{\b}{t}=\f{\a}{\a+\b}\;\;\hskip0.7cm
\ee
such that $A+B=1$.

\begin{figure}[htb]
\centering
\resizebox{0.5\textwidth}{!}{
\includegraphics{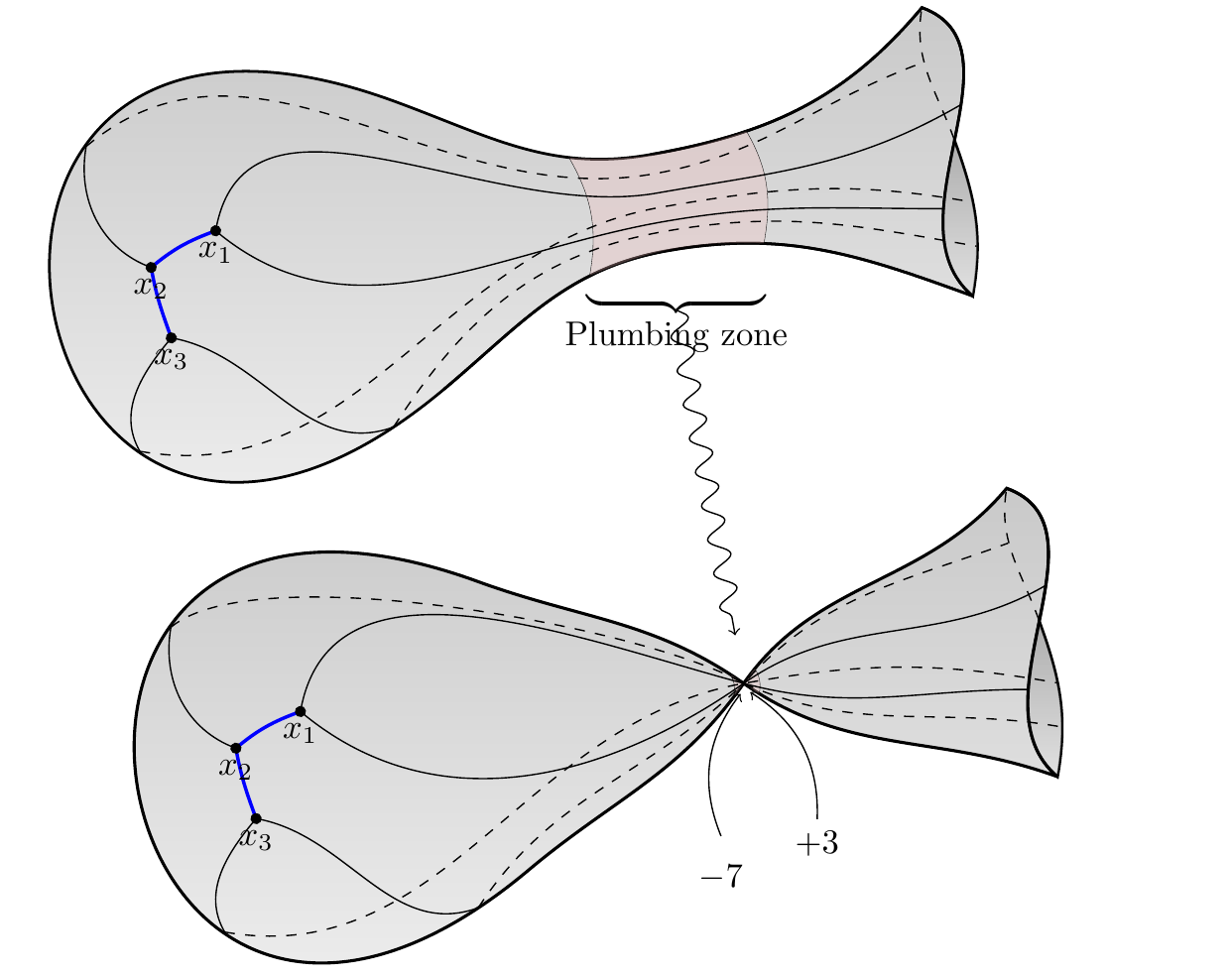}
}
\caption{Separating of the Riemann sphere with three simple zeros and a pole of order 7 of $Q$ in a neighbourhood of $W_1 $}
\la{splitW5}
\end{figure}

Let us now excise from $\Ccal$ a disk $\Dcal_{2t}$ of radius $R=2t$ with centre at $x_1$.

To replace this disk we consider quadratic differential $Q_0$ of the form (\ref{Q0}) such that the absolute value of the periods of the corresponding Abelian differential $v_0$ are given by $2A$ and $2B$ (instead of $2\a$ and $2\b$ as before) so that the lengths of the
finite edges of the corresponding fatgraph equal $A$ and $B$.

In parallel to the construction of the previous section we cut from the Riemann sphere  the 
domain $\Dcal_{2/t}^{A,B}$.
 The  radius of the hole of $\Ccal\setminus \Dcal_{2t}$ is $2t$ while  the radius of the domain
  $\Dcal_{2/t}^{A,B}$  equals $2/t$. Therefore, denoting as before the flat coordinate on $\Ccal$ in a neighbourhood of $x_1$ by $z$, and the flat coordinate near the boundary  of $\Dcal_{2/t}^{A,B}$ by $w$, one has to identify the boundaries via $\frac z w =t^2$.  However, the distinguished local coordinate near $\wt x$ on $\Ccal$ (in terms of the flat  coordinate $z =\int_{\wt x}^x v$ ) is $\zeta=\le(\frac{5}{2} z\ri)^{2/5}$. On the Riemann sphere side, we consider the local coordinate $\xi(x)=\left[\frac {5}2 w(x)  \right]^{-2/5}$ in terms of the global flat coordinate $w(x) = 1+ \int_{x_2}^{x} v_0$ as in  \eqref{flatcoordab}.
  
  Now, consider on $\Ccal\setminus \Dcal_{2t}$ a small annulus $A_{\e}$ defined by $2t<|z(x)|< 2t+\e$ (with coordinate $\zeta$
  on $A_{\e}$).
  On the Riemann sphere side we consider the annulus $A^0_\e$ defined by $2/t< w(x)<2/t+\e/t^2$ (with coordinate $\xi(x)$ 
  on  $A^0_\e$).
 
  Then the identification of the annuli $A_{\e}$ and $A^0_\e$ which is required by the standard plumbing construction 
  is defined by the relation
  \be
  \xi(x) \zeta(y) = t^\frac 4 5 
  \ee
  where $x\in A_\e$, $y\in A^0_\e$.
  
  The "plumbing" point of view is illustrated in Fig. \ref{splitW5}. We notice the following difference of the present framework
  with the standard  plumbing picture. Usually the plumbing parameter is complex and can be used as a local coordinate near the boundary. In our framework we have two real deformation parameters: one of them is the plumbing parameter, 
  while the second is the ratio $\frac B A\in [0,\infty]$ which defines the meromorphic differential on the attached Riemann sphere.

\section{ A local model near Kontsevich's boundary  $W_{-1,-1} $  and space $Q_0^\R([-3]^2)$}
\la{localW11}

The cells of the main stratum  $W_{-1,-1}$  of the Kontsevich boundary   of  the combinatorial model $\Mcal_{g,n}[\pb]$
correspond to fatgraphs which have all vertices of valency 3 except two vertices which have valency 1. The one-valent vertices 	 arise as the result of the normalization of the node.
The corresponding Riemann surface $\Ccal$ is   a (connected or disconnected)  nodal curve  which belongs to Deligne-Mumford boundary of $\Mcal_{g,n}$. The JS differential $Q$ on $\Ccal$ has simple poles at both  points $x_1^0$ and $x_2^0$ obtained by  normalization of the nodal point on the curve  $\Ccal$. The flat metric $|Q(x)|$ has at $x_i^0$ conical points with cone angle $\pi$.

Let us introduce  distinguished local coordinates $\xi_{1,2}$ near $x_{1,2}^0$ such that $Q$ 
in a neighbourhood of $x_i^0$ takes the form
\be
Q(\xi_i)= \frac{(d \xi_i)^2}{\xi_i}\;\;\;\;{\rm as} \;\;\; x\to x_i^0\;.
\ee

The points $x_{1,2}^0$ are branch points on the canonical covering $\Ch$ with distinguished local parameters on $\Ch$ given by $\hat{\xi_i}=\xi_i^{1/2}$; then near $x_{i}^0$ we have $v=2 \d\hat{\xi}_i$. Thus the Abelian differential $v=\sqrt{Q}$ is non-singular and non-vanishing on $\Ch$ at the points
$x_i^0$. 
Flat coordinates near the points $x_i^0$ are given by
$$
z_i(x)=\int_{x_i^0}^x v = 2 \xi_i^{1/2}
$$
such that $Q(x)= (d z_i)^2$ near $x_i^0$.
Conversely, 
$$
\xi_i(x)=\frac{z_i^2}{4}
$$

Denote  lengths of  edges ending at $x_1^0$ and $x_2^0$ by $l_1$ and $l_2$, respectively.
Assume that the   branch cuts ending at $x_1^0$ and $x_2^0$ go along
these edges.

\subsection{Resolution of two one-valent vertices by flat surgery}

To deform the Riemann surface $\Ccal$ into a Riemann surface $\Ccal^{\a,\b}$ which corresponds to JS differential with all simple zeros we fix some $R$ such that $R< \ell_1$ and $R< \ell_2$.
Denote the flat disks around $x_i^0$ defined by inequality $|z_i(x)|<R$ by $\Dcal_{R,i}$ (denote also $\Dcal_R=\Dcal_{R,1}\cup \Dcal_{R,2}$).
 Excising the disks  $\Dcal_{R,i}$ from $\Ccal$  we get an open Riemann surface $\Ccal_R$ 
 with two holes of perimeter $\pi R$ in the metric $|Q|$. 
 The canonical covering $\Ch$ with deleted pre-images of $\Dcal_{R,i}$  turns 
 into 2-sheeted covering $\Ch_R$ of $\Ccal_R$;  $\Ch_R$ is also an open  Riemann surface with 
 two holes.
 
 We are going now to construct a deformation $\Ccal^{\a,\b}$ with two real "small" parameters $\a$ and $\b$ such that the fatgraph corresponding to $\Ccal^{\a,\b}$ has two 3-valent vertices (i.e. two simple zeros
 of the JS differential $Q$) instead of two one-valent vertices $x_1^0$ and $x_2^0$.
 
 Under such deformation two new edges of lengths $\a$ and $\b$ appear while the lengths of edges $\ell_1$ and $\ell_2$ modify as follows:
 \be
 (\ell_1,\ell_2)\to (\ell_1-\frac{1}{2}(\a+\b), \ell_2-\frac{1}{2}(\a+\b))
 \la{l12tran}
 \ee
 so that the  perimeters of both faces containing vertices $x_1^0$ and $x_2^0$ remain unchanged.

To construct $C^{\a,\b}$ we are going to attach  to $\Ccal_R$ an annulus-type domain which we construct as follows.

Consider a point of $\Qcal_0^\R([1]^2,[-3]^2)$ i.e  quadratic differential $Q_0$ on $\CP1$ with two poles of order 3 (which we identify with $x=0$ and $x=\infty$), two simple zeros
and all real periods:
\be
Q_0(x)=\frac{(x-x_1)(x-x_2)}{x^3} (dx)^2
\la{Q011}
\ee
The canonical covering corresponding to the differential $Q_0$ is the elliptic curve $\Ch_0$ with branch points at $x_1$, $x_2$, $0$ and $\infty$ in $x$-plane.
The meromorphic Abelian differential $v_0=\sqrt{Q_0}$:
\be
v_0=\f{\sqrt{(x-x_1)(x-x_2)}}{x^{3/2}} dx
\ee
has second order zeros at $x_1$ and $x_2$ and second order poles at $0$ and $\infty$ on the canonical cover.

Assume that the branch cuts are chosen along horizontal geodesics in the metric $|Q_0|$ connecting $0$ with $x_1$  and $x_2$ with $\infty$.  Approaching $x_1$ along geodesics coming from $x=0$ we can either turn left (this edge we denote by $e_2$)
or right (this edge we denote by $e_1$). The lengths of $e_1$ and $e_2$ will be denoted by $\alpha$ and $\beta$, respectively.

\vskip0.5cm
\begin{figure}[h]
\begin{center}
 \includegraphics{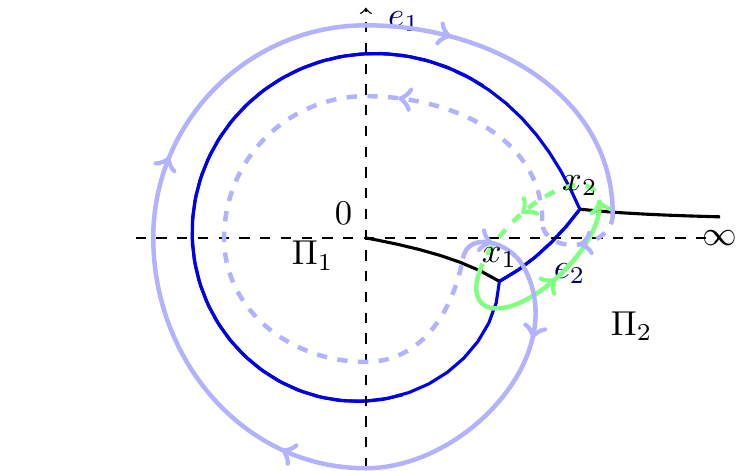}
 \end{center}
\captionof{figure}{Branch-cuts (black) and other horizontal geodesics (blue)  on $\hat \Ccal$. Also shown the edges $e_1$ and $e_2$.}
\label{fig2degmid}
\end{figure}
\vskip0.5cm
The flat coordinate on $\CP1$ and on $\Ch_0$ is given as usual by the Abelian integral of $v_0$; the initial point of integration can be chosen arbitrarily (except $0$ and $\infty$).  
The fatgraph on $\CP1$ corresponding to $Q_0$ has two edges of finite lengths $\a$ and $\b$ which connect points $x_1$ and $x_2$ 
in two different ways, and two edges of infinite length; we enumerate $x_1$ and $x_2$ in such a way that these edges connect $0$ with $x_1$ and $\infty$ with $x_2$ (Fig. \ref{fig2degmid}).\\[3pt]

\noindent \begin{minipage}{0.68\textwidth}
\resizebox{0.99\textwidth}{!}{
\includegraphics{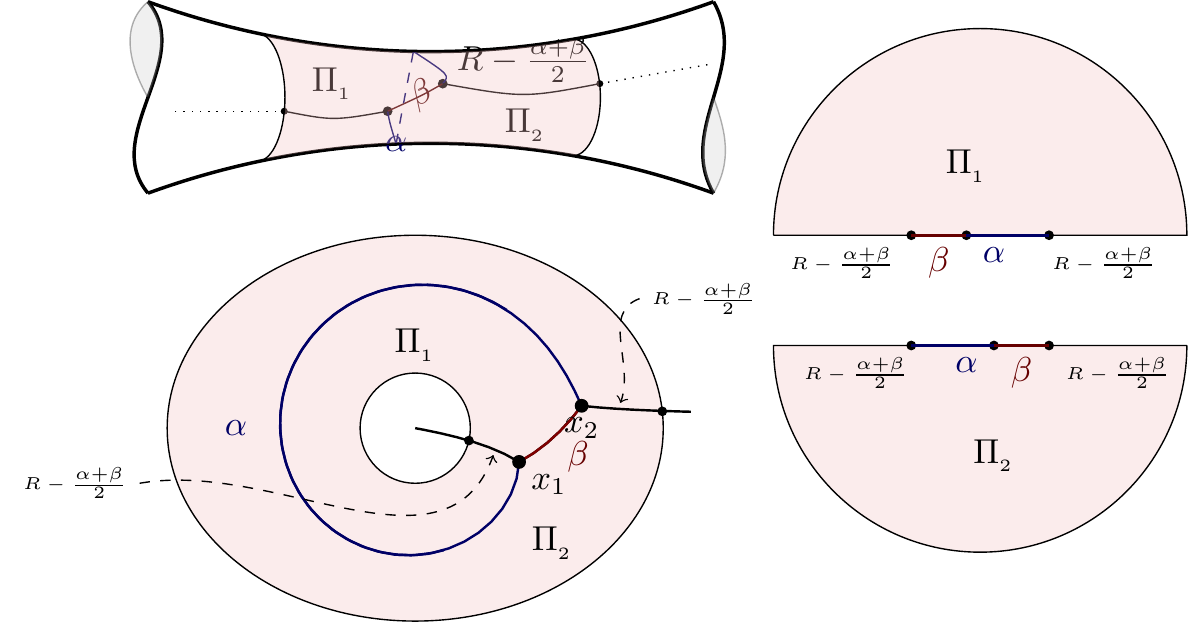}
}
\end{minipage}\,\,\,\,
\begin{minipage}{0.3\textwidth}
\captionof{figure}{The shaded region in $x$-plane is conformally mapped to two half-disks in flat coordinate.  The trajectories in the $x$-plane model are numerically accurate.}
\label{LocalDdeg}
\end{minipage}\\[6pt]
Consider the annular region $\Dcal^{\a,\b}_R$ which consists of two half-disks $\Pi_1$ and $\Pi_2$ of radius $R$ in the flat metric $|Q_0|$,  glued along segments of length $\a$ and $\b$ of their diameters to each other (Fig.\ref{LocalDdeg}). The remaining parts of the diameter 
of each half--disk (having lengths $R-\frac{\a+\b}{2}$) are glued together on each half-disk separately  as shown in Fig.\ref{LocalDdeg}.
The lift   $\widehat{\Dcal}^{\a,\b}_R$ of  $\Dcal^{\a,\b}_R$ to the canonical covering $\Ch_0$ has two outgoing branch cuts coinciding with 
the parts of the diameters of the half-disks of length $R-(\a+\b)/2$. The flat coordinate in $\Pi_1$ 
we  given by $z(x)=\int_{x_1}^x v$ and the flat coordinate in $\Pi_2$ is given by $z(x)=\int_{x_2}^x v$;
then the  boundary of $\Dcal^{\a,\b}_R$ is defined by $|z(x)|=R$.

Perimeters of both boundary components of   $\Dcal^{\a,\b}_R$   equal $\pi R$ in the flat metric $|Q_0|^2$; moreover, both of the boundary components have constant radius $R$ in this metric; the branch cuts go along the real line in flat coordinate. In analogy to the resolution procedure of the 5-valent vertex we "weld"  the region  $\Dcal^{\a,\b}_R$ to two-holed Riemann surface $\Ccal_R$ 
such that the branch cuts on $\Ccal_R$ glue with the outgoing branch cuts on $\Dcal^{\a,\b}_R$.

\begin{figure}[htb]
\begin{center}
\includegraphics{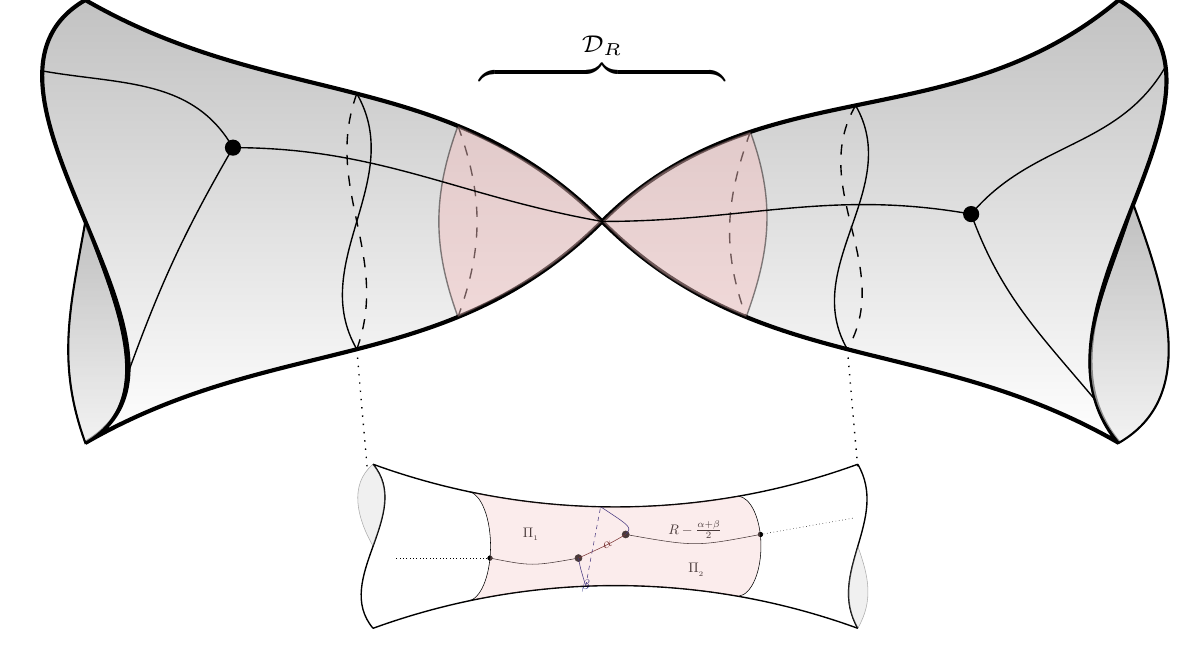}
\end{center}
\captionof{figure}{Resolution of a double point $x_{1,2}^0$ on $\Ccal$ by insertion of an annulus with two simple zeros $x_1$ and $x_2$}
\label{W11plugin}
\end{figure}

The welding produces the smooth Riemann surface $\Ccal^{\a,\b}$ of genus $g$. The fatgraph of $\Ccal^{\a,\b}$ is obtained from the fatgraph of $\Ccal$   by replacing edges of lengths $\ell_1$ and $\ell_2$ on $\Ccal$ 
by edges of lengths $\ell_i-(\a+\b)/2$ (\ref{l12tran}) and introducing two new edges of lengths $\a$ and $\b$. Two one-valent vertices $x_1^0$ and $x_2^0$ on $\Ccal$ are then replaced by two trivalent vertices $x_1$ and $x_2$ on $\Ccal^{\a,\b}$.

\subsection{Plumbing construction}

As well as in the case of the the resolution of a point of $W_1 $,  to study the limit $\a,\b\to 0$ one needs to consider the degeneration of the quadratic 
differential (\ref{Q011}) by taking the limit $x_1, x_2\to 0$.
Alternatively, the "plumbing" construction allows to keep $\Ccal$ and $\Ccal_0$ (together with corresponding JS differentials) fixed while controlling degeneration via
the "plumbing parameter" $t$. The new feature one encounters in the $W_{-1,-1} $ case is the appearance of two plumbing zones instead of one.

To  deform $\Ccal$ into  $\Ccal^{\a,\b}$ (and vice versa, to study the limit $\a,\b\to 0$ of $\Ccal^{\a,\b}$)
we introduce again the plumbing parameter $t=\a+\b$ and define two numbers 
\be
A=\f{\a}{t}=\f{\a}{\a+\b}\;\;\hskip0.7cm
B=\f{\b}{t}=\f{\a}{\a+\b}\;\;\hskip0.7cm
\ee
such that $A+B=1$.

As before, construct the surface $\Ccal_R$ by excising  two disks of radius $R=2t$ from $\Ccal$ (in the flat metric $|Q|$) with centers at $x_{1,2}^0$ so that $R$ is smaller than the lengths $\ell_{1,2}$ of the edges of the 
fatgraph of $\Ccal$ ending at $x^0_{1,2}$.

To replace this disk we consider the quadratic differential $Q_0$ of the form (\ref{Q011}) such that the periods of the corresponding Abelian differential $v_0$ over the cycles $\tilde{a}$ and $\tilde{b}$ shown in Fig. \ref{fig2degmid} are given by $2A$ and $2B$ (instead of $2\a$ and $2\b$ as before) such that the lengths of the
edges connecting $x_1$ and $x_2$ of the corresponding fatgraph on $\CP1$ equal $A$ and $B$.

Similarly to the construction of the previous section we cut from the Riemann sphere  the annular domain $\Dcal_{2/t}^{A,B}$; the lengths of boundaries of $\Dcal_{2/t}^{A,B}$ in the metric $|Q_0|$ equal $2\pi/t$. This time the diameters of the holes of of $\Ccal_{2t}$ differ from diameter of boundaries of
  $\Dcal_{2/t}^{A,B}$  by the factor of $t^2$. Denote as before the flat coordinates on $\Ccal$ in neighbourhoods of $x_{1,2}^0$ by $z_{1,2}$.
  
   The flat coordinates near two boundaries of  $\Dcal_{2/t}^{A,B}$ are denoted by $w_{1,2}$; these coordinates have to be identified with 
  $t^2 z_{1,2}$. However, the proper local coordinates near $x_{1,2}$ on $\Ccal$ are $\zeta_{i}=z_i^2/4$. On the other hand, on the Riemann sphere side, we consider   the local coordinates $\xi_{1,2}(x)=4\left[\int_{x_{1,2}}^x v_0\right]^{-2}$ near the boundaries of   $\Dcal_{2/t}^{A,B}$; the coordinate $z_1$ is used as a flat coordinate 
  in $\Pi_1$ and the coordinate $z_2$  is used as a flat coordinate in $\Pi_2$.

  Now, consider on $\Ccal_{2t}$ two  small annuli ${A_{\e}}_{1,2}$ defined by $2t<z_i(x)< 2t+\e$ (with coordinates $\zeta_i$ on ${A_{\e}}_{1,2}$).
  On the Riemann sphere side we consider the annuli ${A^0_\e}_{1,2}$ defined by $2/t< w_i(x)<2/t+\e/t^2$ (with coordinate $\xi(x)$ on $A^0_\e$).
 
  Then the identification of the annuli ${A_{\e}}_i$ and ${A^0_\e}_i$, required by the standard plumbing construction, is defined by the relation
  \be
  \xi_i(x) \zeta_i(y) = t^{4}
  \ee
  where $x\in {A_\e}_i$, $y\in {A^0_\e}_i$.
  
\begin{figure}[htb]
\begin{center}
\resizebox{0.4\textwidth}{!}{
\includegraphics{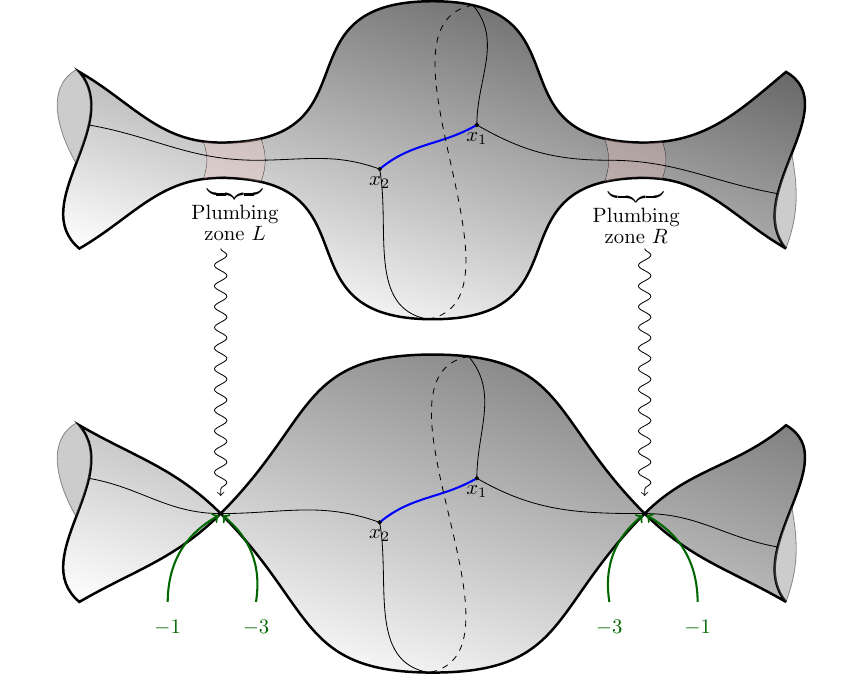}
}
\end{center}
 \captionof{figure}{Plumbing construction interpretation of a resolution of a point  of Kontsevich's boundary $W_{-1,-1}$: a Riemann sphere with two trivalent and two one-valent vertices is glued 
 between two one-valent vertices by introducing two plumbing zones}
 \label{W11-blowup-Fig}
\end{figure}

  The "plumbing" point of view is illustrated in Fig. \ref{W11-blowup-Fig}. Again, the difference with the traditional plumbing construction is that the plumbing parameter $t$ is real, while the second deformation parameter $B/A$ is hidden in the moduli of $\Ccal_0$.

\section{Modular discriminant $\Delta_1 $ on the space  $\Qcal_0^\R(-7)$ }
\la{space5}

 \subsection{Space $\Qcal_0(-7)$}
\label{Q-7}
\paragraph{Stratification.}
 Denote by     $\Qcal_0(-7)$  the complex moduli space of quadratic differentials   $Q$ (we omit the index $0$ on $Q$ from now on) on $\CP1$ with one pole of multiplicity 7: the equivalence is that 
 $Q_1 \sim Q_2$ if there exists a M\"obius transformation $\mu$ such that $\mu^*Q_1=Q_2$.
Therefore, without loss of generality,  we can assume that the pole of $Q$  is  $x=\infty$; by further shift and rescaling of $x$ the differential
 $Q$  can be represented as follows
\be
Q=(x-x_1)(x-x_2)(x-x_3)(dx)^2
\la{Q71int}
\ee
where $x_1+x_2+x_3=0$  and the points are taken up to permutations.   We leave it to the reader to verify that  differentials of the form (\ref{Q71int})  are equivalent iff the sets of zeros 
are given  by $\{x_i\}$ 
and $\{\e{x}_i\}$ where $\e$ is the fifth root of unity.

The space $\Qcal_0(-7)$  has complex  dimension two; it can be stratified according to multiplicities of zeros of the quadratic differential  $Q$ as follows:
\be
\Qcal_0(-7)=\Qcal_0(-7,[1]^3)\cup \Qcal_0(-7,2,1)\cup \Qcal_0(-7,3)\;.
\ee
Here $\Qcal_0(-7,[1]^3)$ is the biggest stratum of  complex dimension 2; it can be identified with 
\be
\Qcal_0(-7,[1]^3)=\{(x_1,x_2,x_3)\in \C^3\;,\;\;\; x_1+x_2+x_3=0\;,\;\; x_i\neq x_j\}/(S_3\times \Z_5)
\la{Q0111}
\ee
Consider the 
canonical cover $\Ch$ defined by the differential $Q$:
\be
y^2=(x-x_1)(x-x_2)(x-x_3)\;.
\la{elliptic}\ee
Introduce on $\Ch$ some basis of canonical cycles $(a,b)$ and define homological coordinates
by integrating the differential $v=y \d x$ over these cycles:
$2A=\int_a v$, $2B=\int_b v$.
The stratum $\Qcal_0(-7,2,1)$ coincides with  the union of   three hyperplanes $x_i=x_j$.
Let  $x_2=x_3$ so that   $x_1=-2x_2$. Then,  
a point of $\Qcal_0(-7,2,1)$ is represented by a differential of the form 
\be
Q=(x+2x_2)(x-x_2)^2(dx)^2\;.
\la{Qwall}\ee
Finally, the stratum  $\Qcal_0(-7,3)$ contains only one point, represented by  the quadratic differential 
\be
Q=x^3 (dx)^2.
\la{QW5}\ee

\paragraph{Coordinatization.}
The periods $(A,B)$ can be used as local coordinates on $\Qcal_0(-7,[1]^3)$ according to the following lemma:
\begin{lemma}
 The Jacobian determinant
of $(A,B)$ with respect to $(x_1,x_2)$ is given by:
\be
\frac{\p(A,B)}{\p(x_1,x_2)}= 2\pi i \Delta_{1}^{1/2}
\la{Jacdet}
\ee
where $\Delta_1 $ is the modular discriminant
\be
\Delta_1 =[(x_1-x_2)(x_1-x_3)(x_2-x_3)]^2
\la{moddis}\ee
\end{lemma}
{\it Proof.}  We have 
\be
\frac{\p(A,B)}{\p(x_1,x_2)}
 = \oint_a  \frac {\pa v}{\pa x_1}  \oint_b \frac {\pa v}{\pa x_2} - \oint_b \frac {\pa v}{\pa x_1}    \oint_a  \frac {\pa v}{\pa x_2}\;.
 \la{rbil}
\ee
To compute this determinant we are going to replace the differential $\p v/\p x_1$  by 
\be
w:=  \frac {\pa v}{\pa x_1}   - \frac { 2x_1+x_2}{x_1+2x_2}\frac    {\pa v}{\pa x_2}  
 = -\frac {(2x_1+x_2)(x_1-x_2)}{2 \sqrt{(x-x_1)(x-x_2)(x-x_3)}} \d x
\ee
which does not change the value of the determinant. Now computation of (\ref{rbil}) reduces to Riemann bilinear relation between 
the holomorphic differential $w$  and $\pa v/\pa x_2$ with the only contribution given by the residue at $x=\infty$.
This residue is easily computed to give (using $t = \frac 1{\sqrt{z}}$ as a local parameter near the infinite point):
$$
\frac {\pa v}{\pa x_2}= \le(\frac {x_1+2x_2}{ t^2} + \mathcal O(1)\ri)\d t\;;
$$
$$
w =\le( (2x_1+x_2)(x_1-x_2) + \mathcal O(t)\ri) \d t\;,
$$
which leads to (\ref{Jacdet}).\QED

 The canonical basis $(a,b)$ is defined up to an $SL(2,\Z)$ transformation. Thus the periods $(2A,2B)$,
 as well as the periods $2\omega_1=\int_a v_0$ and $2\omega_2=\int_b v_0$ 
  of the holomorphic differential $v_0=dx/y$, are also defined up to a $SL(2,\Z)$ action.

The action of $SL(2,\Z)$  on the  ratio $\kappa=B/A$  is the same as   its action on the period
$\sigma=\omega_2/\omega_1$ of the elliptic curve $\Ccal$. In contrast to $\sigma$, which always satisfies 
condition $\Im \sigma>0$ (in particular $\omega_1$ and $\omega_2$ can not be simultaneously real)
no such condition exists for $\kappa$.  

Associating the period $\sigma$  of the corresponding canonical cover to the differential  $Q$  we get a natural fibration of $\Qcal_0(-7,[1]^3)$
over the modular curve  $\Omega$. The fibre can  be identified with $\C^*/\Z_5$ because  the $\C^*$--action  
$(x_1,x_2,x_3) \mapsto (\lambda x_1,\lambda x_2,\lambda x_3), \ \lambda \in \C^*$,  leaves $\sigma$ invariant but produces generically a new point of $\Qcal_0(-7,[1]^3)$, unless $\lambda$ is a fifth root of unity, $\lambda^5=1$.

\subsection{Real  slice $\Qcal_0^\R(-7)$: Boutroux curves}

Elliptic curves corresponding to real periods of the differential $v$ are known as Boutroux curves
\cite{Boutroux, Fokasbook}.

The space $\Qcal_0^\R(-7)$ is the real slice of  $\Qcal_0(-7)$ 
where all the strata  are subsets of the corresponding strata of $\Qcal_0(-7)$ determined by the condition that all periods of $v=y \d x$ 
on the elliptic curve $v^2=Q$ are real. Thus, the space is stratified
as follows:
\be
\Qcal_0^\R(-7)=\underbrace{\Qcal_0^\R(-7,[1]^3)}_{\dim_\R=2}\sqcup \underbrace{\Qcal_0^\R(-7,2,1)}_{\dim_\R=1}\sqcup \underbrace{\Qcal_0^\R(-7,3)}_{\dim_\R=0}
\ee
 Since $\Qcal_0(-7,3)$ consists of only one point,  we have
$\Qcal_0^\R(-7,3)=\Qcal_0(-7,3)$. 

We start from the following lemma describing the configurations of branch points  of
Boutroux curves:

\begin{lemma}[$\Z_5$ and $\R_+$ action]
\label{Z5}
Suppose that $v = \sqrt{ (x- x_1) (x-x_2)(x-x_3)} \d x,\ \  x_1 + x_2 + x_3=0$ has all real periods.
Then so do the differentials
\be
 \sqrt{ (x - {\rm e}^{\frac {2i\pi k}5} \ x_1) (x-{\rm e}^{\frac {2i\pi k}5}  x_2)(x- {\rm e}^{\frac {2i\pi k}5}  x_3)} \d x\;, \ \ \ \  k=1,\dots,4
 \ee
and
\be
\sqrt{ (x - \lambda x_1) (x-\lambda x_2)(x- \lambda x_3)} \d x\;,\ \ \ \  \lambda \in \R _+.
\ee
\end{lemma}
{\bf Proof.}
Under the map $x_i\mapsto \lambda {\rm e}^{\frac {2i\pi k}5} x_i$, $i=1,2,3$   the periods transform as $P\mapsto \lambda^\frac 52 {\rm e}^{ {i\pi k}}  P $.  Hence under this map the  periods are $\R$-scaled  but  remain real. 
\QED

\begin{lemma}\la{LemBout}
{\bf [1]} If  $Q\in \Qcal_0^\R(-7,[1]^3)$ then  $x_1$, $x_2$ and $x_3$ 
do not lie on the same line through the origin in the $x$-plane.\\{}
{\bf [2]} If  $Q\in \Qcal_0^\R(-7,2,1)$ 
then ${\rm arg}\; x_{1,2}= \le\{\frac{2i\pi k}5,\ k=0,\dots 4\ri\}$ and  ${\rm arg}\; x_{3}={\rm arg} \;x_{1,2}+\pi$.
\end{lemma}
{\it Proof.}
{\bf [1]} Let us show that the assumption that all $x_i$ lie on the same line ${\rm e}^{i\theta} \R$ is 
incompatible with the reality of all periods of $v$.
Using an $\R_+$-rescaling one can assume without loss of generality that  $x_1 = {\rm e}^{i\theta}$, $x_2 = \lambda x_1$ and therefore $x_3 = -(1+\lambda) x_1$ with some $\lambda>1$. 
Then the periods of $v$ are given by  the integrals
\be
\oint_\gamma v = {\rm e}^{\frac 52 i\theta} \oint_\gamma  \sqrt{(\tilde{x}-1)(\tilde{x}-\lambda)(\tilde{x}+1+\lambda)} \d \tilde{x}\;,
\la{intv10}
\ee
where $\gamma$ is an arbitrary cycle going around the branch points $1$, $\lambda$ and $-1-\lambda$ in the
$\tilde{x}$-plane.
The integral around  the  cut $[1,\lambda]$ in (\ref{intv10}) is imaginary and the integral around  
the interval $[-1-\lambda,1]$ is real. Thus there is no such value of $\theta$ that both of these periods are simultaneously real unless some of the branch points coincide.
When some of the branch points coincide, say,  $\lambda=1$, then the reality of the expression
${\rm e}^{\frac 52 i\theta} \int_{-2}^1 (\tilde{x}-1) \sqrt{(\tilde{x}+2)} \d \tilde{x}$
implies  $\frac 5 2\theta = 0 \,{\rm mod} \, \pi $, which proves part {\bf [2]}.
\QED

\begin{lemma}
\label{Lemmaemer}
Let $Q \in \mathcal Q_0^\R(-7,2,1)$ so that the double root lies on one of the rays $ {\rm e}^{\frac {2i\pi k}5 }\R_+$. Then, for any  infinitesimal deformation of $Q$ preserving the
 Boutroux condition, the double root splits into two roots along one of the two directions forming an angle $\pm \frac \pi 4$ with the ray $ {\rm e}^{\frac {2i\pi k}5 }\R_+$; the central root emerges along the direction $\pm \frac 54\pi$, with the angle measured  from the ray $ {\rm e}^{\frac {2i\pi k}5 }\R_+$ with the positive orientation.
\end{lemma}
{\bf Proof}.
Denote by $x_{2,3}$ the two zeroes emerging from the double zero and denote by $A$ the period of $v$ along a small cycle encircling  them. Using the $\Z_5$ action (Lemma \ref{Z5}) we can assume without loss of generality that if $x_2=x_3$ then ${\rm arg} \,x_{2,3}=2\pi/5$ and ${\rm arg} \,x_{1}=7\pi/5$. 
To compute the limit of ${\rm arg}(x_3-x_2)$ as $x_3\to x_2$ we proceed as follows. Let $x_3=x_2+\delta$,
and $x_1=-2x_2 -\delta$ where $x_2(\delta)= x_0 + \Ocal(\delta)$ for some $x_0=|x_0| e^{2\pi i/5}$.
Then, up to higher order powers of  $\delta$,
$$
A= \pm \int_{x_2}^{x_2+\delta}[(x-x_2)(x-x_2-\delta)(x+2 x_2 +\delta)]^{1/2}\d x
$$
$$
\sim \pm (3 x_0)^{1/2}
\int_{x_0}^{x_0+\delta}[(x-x_0)(x-x_0-\delta)]^{1/2}\d x=\pm (-3 x_0)^{1/2} \delta^2
$$
Therefore, since $A\in \R$ and $ \arg(x_0)= 2\pi/5$, we must  have (as $\delta \to 0$)
$$
2 {\rm arg}\,\delta +\f{\pi}{2} +\frac{\pi}{5} \to 0 \;\;\; ({\rm mod}\;\; \pi)
$$
and, hence ${\rm arg}\,\delta\to \f{2\pi}{5}\pm\f{\pi}{4}$ as $\delta\to 0$. This proves the first assertion. 

To prove the second assertion we  need to analyze the critical trajectories emerging from $x_{2,3}$ for small $\delta$. To this end we define  
\be
\z ={\rm e}^{-\frac {2i\pi}5}  \frac{(x-|x_0|)}{|\delta|}  , \qquad
q(\z;|\delta|) = |\delta|^{-4} Q(x)  
\ee
so that $x_2, x_3$ are mapped to two points at distance $1$ from the origin of the $\z$--plane, whereas $x_1$ is mapped to a point at distance $\mathcal O(\delta^{-1})$ from the origin (see Fig. \ref{figdege}). 
Then, a straightforward computation yields
\bea
q(\z;|\delta|) 
 = \le (|\delta| \z  + 3r \ri) \le(\z^2  - {\rm e}^{2i\arg\delta- \frac {4i\pi}5}\ri) \d \z^2\;.
\eea
Note that $q$ is also  quadratic differential with real periods and its trajectory structure has the same topology as the trajectory structure of  $Q$ (up to an affine transformation). The ray $ {\rm e}^{\frac {2i\pi }5 }\R_+$ is mapped to the positive real axis of the $\z$--plane. 
Now, as $ |\delta|\to 0$ (i.e. $A\to 0 $) we obtain (uniformly over compact sets in the $\z$--plane)
\be
\label{locscal}
\lim _{A\to 0_+} q(\z;|\delta|) = 3 r (\z^2 \pm i)\d \z^2  \;.
\ee
Consider the case $-i$ in the above expression, the $+i$ case being treated similarly.
A simple analysis of the structure of geodesics in the metric given by the modulus of the quadratic differential  (\ref{locscal}) shows  that the trajectory that extends to $\Re \z = -\infty$ issues from ${\rm e}^{\frac {5i\pi}4}$ (see Fig. \ref{figdege}).
Thus the central  zero $x_2$ is the one that is mapped to ${\rm e}^{\frac {5i\pi}4}$. Similarly for the case $-$. 
\QED

The reality of the periods of $v$ on $\Ch$ allows to connect $x_1$, $x_2$ and $x_3$ by two horizontal trajectories 
of the differential $v$. Then one of the zeros (say, $x_2$) is connected to $x_1$ and $x_3$ and will be referred to henceforth as the  {\it central
zero}. Two other zeros $x_3$ and $x_1$ will be labeled according to values of their arguments: counterclockwise with respect to the origin we assume that ${\rm arg} x_2 < {\rm arg} x_3< {\rm arg} x_1$. 

All periods of $v$ are real. In addition, it turns out to be convenient to work with positive numbers; thus we introduce 
an additional modulus in the real setting and define 
\be
A=\Big|\int_{x_2}^{x_3}v\Big|\;;\hskip0.7cm B=\Big|\int_{x_2}^{x_1}v\Big|
\la{realper}
\ee
where the integration contours go along the horizontal trajectories of $v$ and correspond to $a/2$ and $b/2$
in homologies. Therefore, $A$ and $B$ are the lengths of horizontal geodesics connecting $x_2$ with $x_3$ and $x_1$, respectively.

The flat coordinate $z(x)$ is defined on $\Ch$ by integration of $v$ with initial point $x_2$:
\be
z(x)=\int_{x_2}^x v
\la{flatcoord}
\ee
with the branch cuts chosen as   the union of trajectories $[x_2,\infty)$,  $[x_1,x_2]$ and $[x_2,x_3]$.

The following theorem is an analog of Strebel's theorem which allows to reproduce a Riemann surface with punctures knowing the
lengths of edges and topology of the corresponding fatgraph.

\begin{theorem}
\label{thm123}
{\bf [1]} The ordered pair of lengths $(A,B)$  (\ref{realper}) defines a one-to-one map of  the space 
$ \Qcal_0^\R(-7,[1]^3)$ to $\R_+^2$.
\\{}{\bf [2]}
The diagonal $A=B$ corresponds to   Boutroux curves with an antiholomorphic involution (the {\rm Boutroux-Krichever locus}).
By the  $\Z_5$ action of Lemma \ref{Z5},  
the central zero $x_2$ of a Boutroux-Krichever curve can always be assumed to be real and positive while the other zeros form a conjugated pair: $x_1=\bar{x_3}$.
\\{} {\bf [3]}  
The boundaries $A=0$ and $B=0$, being identified, correspond to points of the space $ \Qcal_0^\R(-7,2,1)$.
The remaining non-vanishing length parametrizes $ \Qcal_0^\R(-7,2,1)$.
\end{theorem}
{\it Proof.}   {\bf [1]} The Jacobian of coordinates $(A,B)$ on $\wh \Qcal_0^\R (-7,[1]^3)$ defined by \eqref{realper} never vanishes according to (\ref{Jacdet});
however this is not sufficient to prove that a given function of two variables is globally one-to-one.
To prove   the  invertibility  we need to restore the triple $(x_1,x_2,x_3)$ 
uniquely (up to the $\Z_5$ action) from the knowledge of the periods $A$ and $B$ (\ref{realper}).
Knowing the lengths $A$ and $B$ of the critical graph we construct a polyhedral surface (i.e. surface with flat metric and conical singularities) by gluing 5 half-planes as shown in Fig.\ref{localmodelW5}. This  uniquely defines the conformal structure i.e 
$J$-invariant of $\Ch$. Moreover, since the zero $x_2$ is assumed to be "central" and the zeros $x_1$ and $x_3$ are labeled we
can choose a distinguished pair of canonical $a$- and $b$-cycles by assuming that $a$ encircles the critical trajectory connecting $x_2$ and $x_3$. Moreover, we fix the direction of the $a$-cycle such that $A=\f{1}{2}\int_a v$.  The  $b$-cycle then encircles the critical trajectory connecting $x_2$ and $x_1$, with an appropriate orientation. This determines the period $\sigma$ of
$\Ch$ corresponding to such distinguished Torelli marking. In turn, this allows to express the cross-ratio of the points $x_i$ in terms of corresponding theta-constants by the use of Thom\ae\ formulas. We call this ratio $t$:
\be
\f{x_2-x_1}{x_3-x_1}=t
\la{deft}
\ee
Introducing the new variable $\xt=\frac{x-x_1}{x_3-x_1}$ we can rewrite  the first of the equations (\ref{realper}) as follows:
\be
(x_3-x_1)^{5/2}\int_t^1\sqrt{\xt(\xt-1)(\xt-t)}d\xt = A
\ee
which implies the relation  between $x_1$ and $x_3$ 
\be
x_3-x_1=\left[\frac{A}{\int_t^1\sqrt{\xt(\xt-1)(\xt-t)}d\xt}\right]^{2/5}
\la{x31}
\ee
where the right-hand side is known up to a power  of $e^{4\pi i/5}$.
This allows us to conclude that 
the system of three linear equations (\ref{deft}), (\ref{x31}) and $x_1+x_2+x_3=0$ define all  $x_i$ uniquely up to the $\Z_5$ action,
which proves that the period map provides an isomorphism between the space $ \wh{\Qcal}_0^\R(-7,[1]^3)$ and the first quadrant
in the $(A,B)$-plane.

{\bf [2]} Let  $A=B$. Then the elliptic curve $C$ admits an anti-holomorphic involution $\tau$ which acts as $z\to \bar{z}$ in the  flat coordinate $z$
such that the differential $v=\d z$ satisfies the relation $\overline{v(x^\tau)}=v(x)$.
The  involution $\tau$ must act in the $x$-plane as reflection with respect to some line passing through the origin (since $x_1+x_2+x_3=0$  and also $x=\infty$ must be invariant under $\tau$); thus $\tau$ must interchange $x_1$ and $x_3$ while leaving $x_2$ invariant (all $x_i$ can not be invariant under this involution since they don't lie on the same line as per Lemma \ref{LemBout}, item [1]). 

Therefore, there exists a real $\theta$ such that $x_i= e^{i\theta} x_i'$ while
$\overline{x_1'}=x_3'$ and $x_2'\in \R$.  Making the change of variable $x=e^{i\theta}x'$ we see that $\tau$ acts on $x'$ as $x'^{\tau}=\bar{x}'$;
thus $\tau(x)=e^{-2i\theta} \bar{x}$. 
Then 
$$
Q(x)=(x-e^{i\theta} x_1')(x-e^{i\theta} \bar{x}_1')(x-e^{i\theta} x_2')(dx)^2
$$
and
$$
\overline{Q(x^\tau)}=e^{10i\theta} Q(x)
$$
Since we have assumed that not only $Q$, but also $v=\d z=\sqrt{Q}$ satisfies the relation $\overline{v(x^\tau)}=v(x)$ then $5\theta$ must be a multiple of $2\pi$. Moreover, since we assumed that for non-degenerate curves $ 0 \leq {\rm arg}\,x_2<2\pi /5$ we conclude that $x_2\in \R_+$ 
and $x_3=\bar{x}_1$.

{\bf [3]} Consider a point of $ \Qcal_0^\R(-7,2,1)$ where, say,  $x_2=x_3 = r{\rm e}^{\frac {2i\pi}5}$ and 
$ x_1 =-2 r {\rm e}^{\frac {2i\pi}5} $. Then the remaining non-vanishing period $B$ is explicitly computable:
\bea
B= \Big|\int_{x_2}^{x_1} (x-r{\rm e}^{\frac {2i\pi}5}) \sqrt{ x-2 r {\rm e}^{-\frac {3i\pi}5}} \d x\Big| = \Big|\int_{-2r}^r (r-z) \sqrt{z+2r} \d z\Big| = \frac {12\sqrt{3}}5 r^\frac 5 2\;.
\eea 
This shows that the period $B\in \R_+$ is a global coordinate on $ \Qcal_0^\R(-7,2,1)$.
\QED

\begin{figure}
\begin{center}
\includegraphics{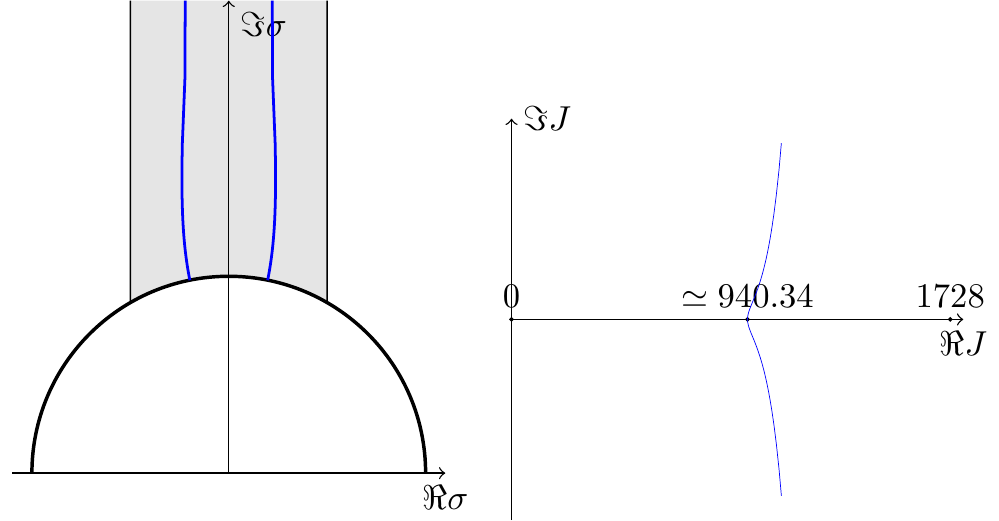}
\end{center}
\caption{Left: The space of Boutroux curves $\Qcal_0^\R(-7)$ is  fibered over the shown curve $\Rcal_1$ in the moduli space of elliptic curves with fiber $\R_+$. Right: 
 the set $\Rcal_1$  in the plane of $J$-invariant. The point $J=1728$ corresponds to $\sigma=i$ while  $\sigma=e^{2\pi i/3}$ corresponds to $J=0$.  Note that $1728$ is the year the English astronomer James Bradley gave the first (and amazingly precise!) estimate of the speed of light.}
\la{sigma71int}
\end{figure}


The $J$-invariant of Boutroux curves traces  the curve $\Rcal_1$, computed numerically and shown in Fig.\ref{sigma71int} (right).
The point of intersection of  $\Rcal_1$  with
the real $J$--line  corresponds to the only curve which is simultaneously real (i.e. it admits an anti-holomorphic involution) and  Boutroux. For this curve (known as of Boutroux-Krichever \cite{Wiegman})
the $J$-invariant  is approximately given by $940.34$. The values of $x_i$ corresponding to this curve are given by
\be
x_1=-1- 6.492 i \;,\hskip0.7cm x_2= 2\;,   \hskip0.7cm x_3=-1+6.492 i
\ee
(or any multiple of these values with a positive real constant).
The set $\Rcal_1$ in the plane of the period $\sigma$ of the curve $\C$ is shown in Figure \ref{sigma71int} (left).

Let us discuss now which configurations of zeros $(x_1,x_2,x_3)$ correspond to a given pair of periods
 $(A,B)$. 
 The following Proposition provides the converse to certain statements of Lemma \ref{LemBout} and Thm. \ref{thm123}.

 \begin{proposition}
\label{ThmW5}
Suppose that $v = \sqrt{ (x-x_1) (x-x_2)(x-x_3)} \d x,\ \  x_1 + x_2 + x_3=0$ has all real periods. 
Then  one of the roots $x_i$  belongs to one of the rays ${\rm e}^{\frac {2i\pi k }5} \R_+$  iff   either
\begin{enumerate}
\item  $x_i$ is a double root (i.e. the curve is degenerate) or;
\item  the curve is nondegenerate and admits an anti-holomorphic involution exchanging the other two roots, i.e.  $\wh{ \mathcal C}$   is the Boutroux-Krichever curve of Thm. \ref{thm123} (point [2]). In this case $x_i$ is the central zero, which we denote by $x_2$.
\end{enumerate}
\end{proposition}
{\bf Proof.}
{\it  Part {"$ \Leftarrow$"}.} 
Suppose that the curve is Boutroux-Krichever. Then the  central zero belongs to one of the rays  ${\rm e}^{\frac {2i\pi k }5} \R_+$ by Thm. \ref{thm123}, item [2]. Suppose  the curve  is degenerate. Write the differential as $v = (x-x_2)\sqrt{x+2x_2} \d x$. Then a direct computation gives 
$\int_{x_2}^{-2x_2} v =  \frac {12\sqrt{3}}5 x_2^\frac 52$, and this is real iff $x_2$ belongs to one of the rays.

\paragraph{\it Part {"$ \Rightarrow$"}.} 
Using the $\Z_5$ action and $\R_+$ - scaling  of Lemma \ref{Z5}, we can assume that the root in question is $x_j=2$. Then the other two roots (since the sum is zero)  must be of the form $-1\pm \rho$, for some $\rho\in \C$. Let us write
\be
v= \sqrt{(x-2)\le(x+ 1   +\rho \ri)\le(x+1  - \rho\ri)}\ .
\ee
The curve is  degenerate when  $\rho = 0, \pm 3$.  The case $\rho=0$ is excluded since the periods are not real.  In both cases $\rho= \pm 3$ the periods are real.
 
Let now the curve be non-degenerate:  it is our goal to show that $\rho \in i\R$, which would also imply that $x_j =x_2=2$ is the central root. 

\begin{figure}
\begin{center}
\includegraphics{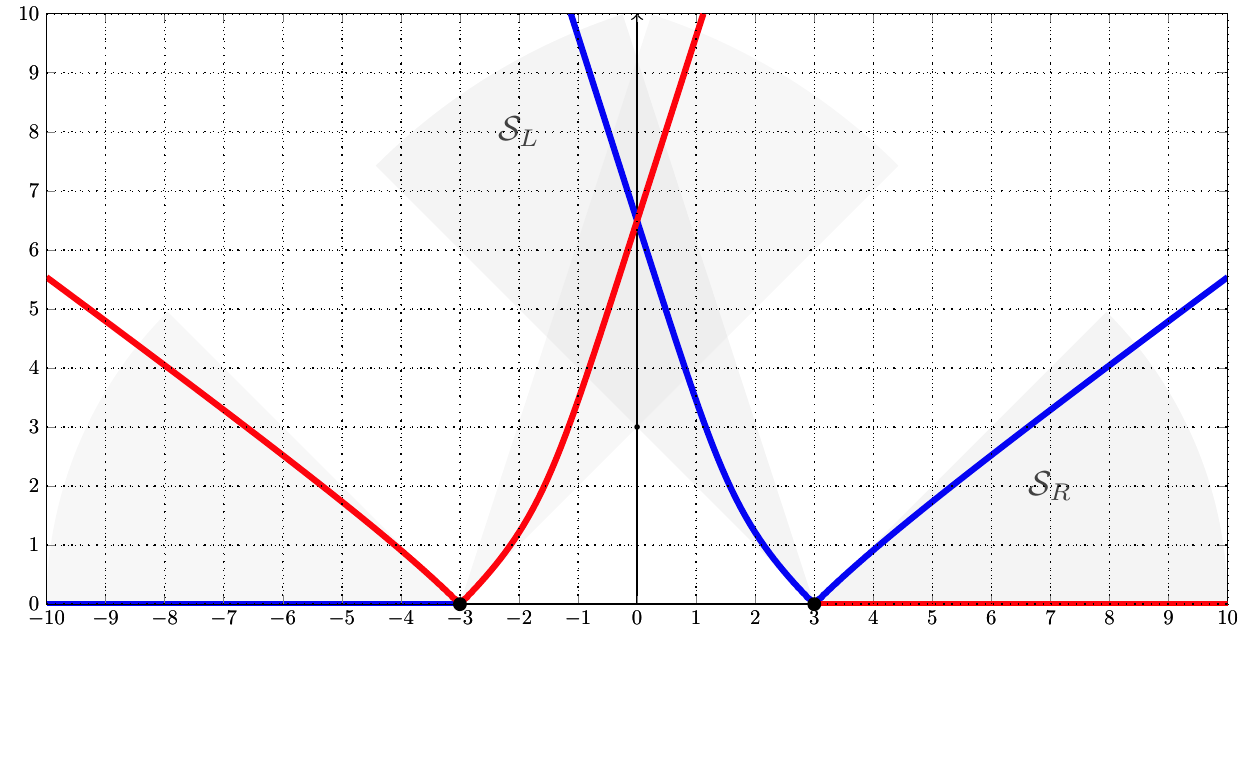}
\end{center}
\vspace{-2cm}
\caption{The zero level sets of $\Im(A_+(\rho))$ (blue) and $\Im (A_-(\rho))$ (red), plotted numerically.}
\label{figIm}
\end{figure}

\begin{figure}[t]
\begin{center}
\resizebox{0.9\textwidth}{!}{
\includegraphics{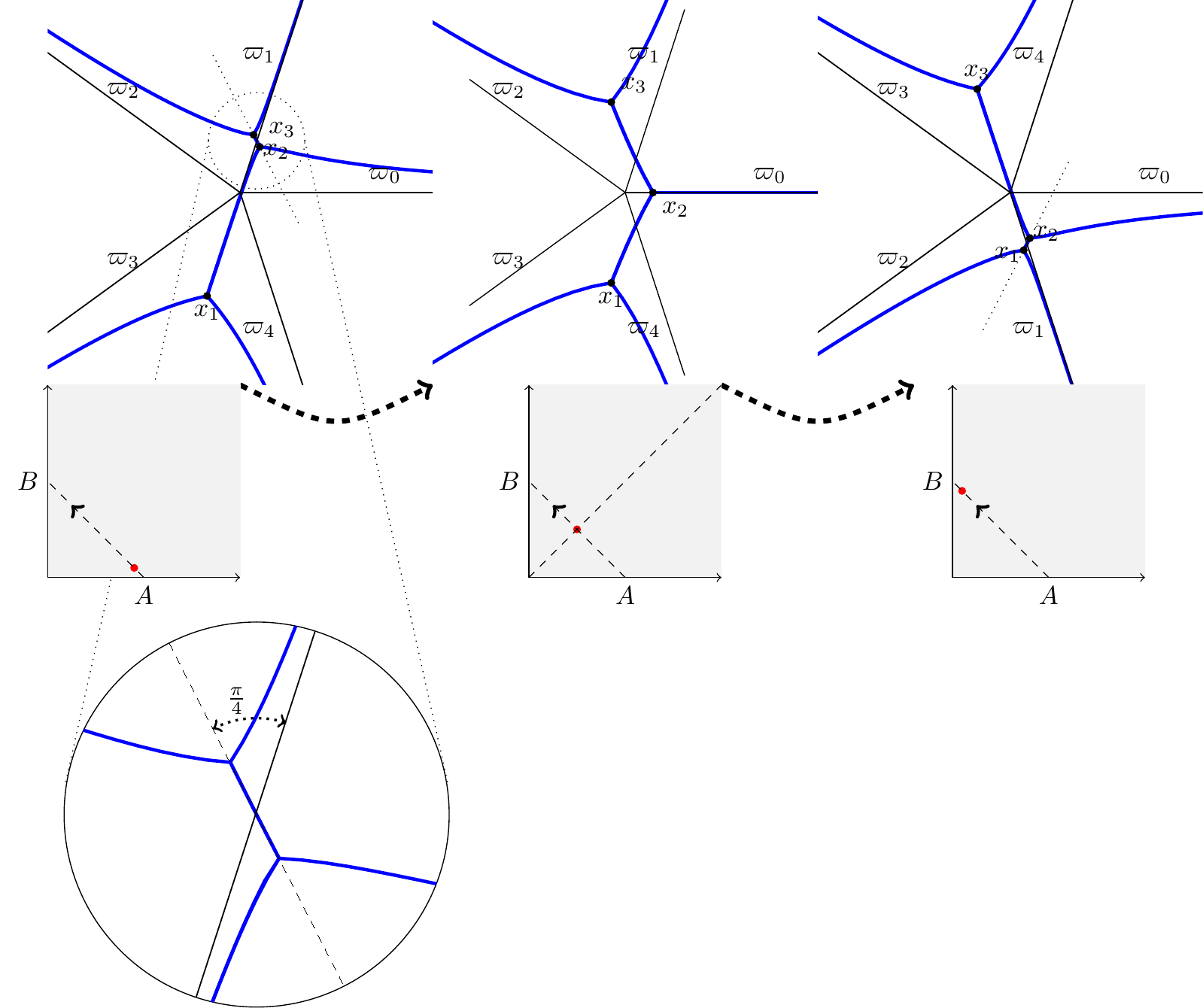}
}
\end{center} 
\captionof{figure}{ The motion of the  roots $\{x_i\}$ as the periods $(A,B)$ traverse a cell from $(0,1)$ to $(1,0)$ along
$A+B=1$.   The zoom-in inset shows the angle of approach of the roots $x_1$ and $x_3$ 
to the  ray ${\rm arg}\,x=2\pi/5$ as $(A,B)\to (0,1)$.}
\label{figdege}
\end{figure}

To find out for which values of $\rho\neq 0,\pm 3$  both periods are real consider the integrals
\be
A_\pm(\rho) =\int_{2}^{-1 \pm  \rho} v \d x
\ee
such that that $A_-(\rho) = A_+(-\rho)$. 
To prove the theorem it is sufficient to show that there is  only one  solution of the system
\be
\Im A_+ (\rho)=\Im A_- (\rho)=0
\la{AA}\ee
in the upper half plane and this solution is  located on $\Re \rho=0$. 
The function  $A_+(\rho)$ is analytic on the domain $\mathcal D_+ = \C \setminus (-\infty,0]$ while  $A_-(\rho)$ is analytic on the domain $\mathcal D_- = \C \setminus [0,\infty)$  and both satisfy the Schwartz symmetry $\ov{A_\pm(\ov \rho)} = A_\pm(\rho)$.
The logic of the proof is as follows:
\begin{enumerate}
\item  Consider the curve (which may consist of several components) $\Phi=0$, with $\Phi=\Im A_+(\rho)$; the solutions of the system (\ref{AA})  lie at the intersection of the two curves $\{\Phi=0\}$ and $-\{\Phi=0\}$. Taking into account  the Schwartz symmetry of the functions $A_\pm(\rho)$ it is sufficient to analyze the upper half plane.
\item  Restricting ourselves  to the upper half-plane we  show that the set $\Phi=0$ is contained in the union of  two sectors, ${\mathcal S_R}$ and ${\mathcal S_L}$ (shown in grey in  Fig. \ref{figIm}) which are defined as follows:
\be
\{\rho:\ \Phi(\rho)=0\} \subset\underbrace{ \le\{\arg(\rho-3)\in \le[\frac {3\pi}5, \frac {3\pi}4\ri]\ri\}}_{\mathcal S_L} \cup \underbrace{\le \{\arg(\rho-3)\in \le[0, \frac \pi 4\ri]\ri\}}_{\mathcal S_R}
\label{SLSR}
\ee
and each of the two sectors contains exactly one smooth branch of the set $\Phi=0$ starting from $\rho=3$.
\item Since $\mathcal S_L \cap (-\ov {\mathcal S_L}) =\emptyset = \mathcal S_R \cap (-\ov{\mathcal S_L}) = (-\ov{\mathcal S_R})\cap \mathcal S_L$, all solutions of (\ref{AA}) must lie in $\mathcal S_L\cap (-\ov{\mathcal S_L})$ (the diamond shaped region in Fig. \ref{figIm}).
\item Show that there is only one point of intersection of $\Phi=0$ with the positive imaginary axis, and this point lies within the diamond region.

\item Finally, prove that within $\mathcal S_L$ the branch of the  $\Phi=0$ is a function of $y = \Im \rho$ by showing that no tangent of the foliation $\Phi=$constant can be horizontal within $\mathcal S_L$.   Together with the previous item this  immediately  implies the uniqueness of the  solution of (\ref{AA}).  
\end{enumerate}

Let us now discuss these items in more detail.

The item {\bf (1.)} is self-evident. 
Let us look at item {\bf (2.)} 
Parametrizing the integral by $x = 2 + (\rho-3)\frac {1-s}{1+s}$, $s\in [0,1]$  we obtain the  integral 
\bea
A_+(\rho)
 \label{A+}
 =i(\rho-3)^2 \int_0^1 \sqrt{\rho+3s}\frac{ \sqrt{s(1-s^2)}}{(1+s)^4}\d s\\
 =
i (\rho-3)^\frac 5 2    \int_0^1 \sqrt{1+3\frac{s+1}{\rho-3}}\frac{ \sqrt{s(1-s^2)}}{(1+s)^4}\d s\;.
 \label{A++}
\eea

From representation \eqref{A++} we see that for $|\rho|\to +\infty$ the four branches have slopes $\pm \frac {\pi}5, \pm \frac {3\pi}5$. 
Moreover, there is a fifth branch consisting of the ray $\rho \in (-\infty, -3]$ (the integral is purely imaginary and hence the expression \eqref{A++} is real)  but this fifth branch does not play a role in this discussion.
	
The representation  \eqref{A++} is also useful in analysing the argument of $A_+$ in the finite region; to this end, let us consider this argument for $\Im \rho>0$;
\be
\label{arga+}
\arg(A_+(\rho))  = \frac \pi 2   + \frac 5 2 \arg(\rho-3) + \arg
\le(  \int_0^1 \sqrt{1+3\frac{s+1}{\rho-3}}\frac{ \sqrt{s(1-s^2)}}{(1+s)^4}\d s\ri)\;.
\ee
It is our goal  to see where the argument of $A_+$ can be equal to $0$ or $\pi$; if the last term in \eqref{arga+} were absent then  the set $\Phi=0$ would simply consist of the rays $\arg(\rho-3) \in\le\{  \frac {3\pi} 5, \frac {\pi}5\ri\}$.   
Observing that integrals with positive integration measures preserve cones, we estimate the last term in \eqref{arga+} to belong to the interval
$(-\frac {\arg(\rho-3)}2, 0)$.
From \eqref{arga+}  and the estimate of the argument of the integral we have 
\be
\label{bounds}
\frac \pi 2 +2\arg(\rho-3)< \arg A_+(\rho) < \frac \pi 2 + \frac 5 2 \arg(\rho-3).
\ee
To obtain the sectors $\mathcal S_{L,R}$ one has to check that in the complement of those sectors the argument of $A_+(\rho)$ cannot be equal to $0$ or $\pi$. For example in the sector $\arg(\rho-3)\in \le[\frac \pi  4, \frac {3\pi}5 \ri]$ lying between $\mathcal S_L$ and $\mathcal S_R$ we obtain from \eqref{bounds} that  $\pi < \arg A_+< 2\pi$ and hence $A_+$ cannot be real there. Similarly, in the sector $\arg(\rho-3)\in \le[\frac {3\pi}4, \pi\ri]$ we have $2\pi < \arg A_+< 3\pi$, and thus $A_+$ cannot be real there as well.
To summarize:
\be
\label{conc1}
\hbox{the set $\Phi=0$ i.e. $\arg A_+(\rho)\in\{0,\pi\}$ lies in the  union of two sectors, $\mathcal S_L$ and $ \mathcal S_R$ \eqref{SLSR}.} 
\ee
(see Fig. \ref{figIm}).
Moreover, in each of the two sectors there is exactly one branch of the curve $\Phi=0$ that extends from $\rho=3$ to infinity; this statement holds since the level sets of harmonic functions cannot form bounded curves unless they  surround a singularity;  meanwhile $A_+(\rho)$ is harmonic in the whole upper half plane. Moreover, expression \eqref{A+} near $\rho=3$ shows that there are  four branches of $\Phi =0$ issuing from $\rho=3$ with slopes $\pm \frac \pi 4  \pm \pi$.
This proves point {\bf(2.)} and hence also {\bf (3.)}.

Therefore, all  solutions to the system \eqref{AA} lie in the diamond--shaped region of intersection $ (-\ov{\mathcal S_L} ) \cap \mathcal S_L$ (see Fig. \ref{figIm}). This region lies in the half-plane   $\Im \rho\geq 3$. 

{\bf (4.)} Let us now show that the set $\Phi=0$ for $\Im \rho>3$  intersects the imaginary $\rho$--axis $i\R_+$ exactly once at some $\rho = i\rho_\star$. This means that for $\Im \rho<\rho_\star$ the set $\Phi=0$ lies entirely in the right half-plane, while the branch in $\mathcal S_L$ for $\Im \rho>\rho_\star$ lies entirely in the left half-plane. In other words the branch of interest does not ``zig-zag'' intersecting  $i\R_+$ several times. 

To this end, let us write out $A_+(iy)$ more explicitly
\be
\Phi(iy) = \Im A_+(iy)\ ,\ \ \ 
A_+(iy) = -i (iy-3)^2 \int_0^1 \sqrt{3s+iy} \d\mu(s) \, , \ \  \  \d\mu(s) = \frac { \sqrt{s(1-s^2)}}{(1+s)^4} \d s\; .        
\ee
Then
\be
\Phi(3i ) =\Im \le( -9i  (i-1)^2 \int_0^1 \sqrt{3s + 3i}\, \d\mu(s)\ri) 
= -18\int_0^1 \Im \le(\sqrt{3s + 3i}\ri)\, \d\mu(s)
 < 0 
\ee
(the approximate value is $\Phi(3i)\simeq -2.2438$).
Moreover, for $y\to+\infty$ we have  $\Phi(iy) \sim   y^\frac 52 \frac {\sqrt{2}}{2} \int_0^1 \d \mu(s)\to +\infty$   so that there is at least one solution $\Phi(iy)=0$ for $y\in [3,\infty)$.

To show uniqueness of this solution suppose that there are several solutions of (\ref{AA})  lying on the imaginary axis; at each of these points both $A_\pm$ are real. Then for each of these points we would have $A=B$ and each would give some  Krichever--Boutroux curve. But we know from Thm. \ref{thm123}$_{[2]}$  that the Krichever--Boutroux curve is unique (up to scaling, which has  already been used to set the central root at $2$) which leads to a contradiction.

{\bf (5.)}
 Let us study the slopes of the tangent vectors to the level curves  and show that their tangent can never be horizontal; in other words here we  show that within $\mathcal S_L$ the argument of $ {\frac {\d  A_+}{\d \rho}} $  never equals $0$ or $\pi$. 
To this end we analyze the argument of the  expression 
\bea
\frac {\d A_+(\rho)} {\d \rho} = \oint \frac { \d v}{\d \rho} \d z
=  i\rho (\rho-3) \int_0^1 \frac {t \d t}{ \sqrt{t (1-t^2) \le(\rho  + 3\frac{1-t}{t+1} \ri)}}
=
\nonumber \\
\label{dAplus}
=
- i\rho(\rho-3) \int_0^1 \frac 1 {\sqrt{\rho+ 3s }} \frac {(1-s) \d s}{ (s+1)\sqrt{s(1-s^2)}} 
\eea
in the sector $S_L= \le\{\arg(\rho-3) \in \le[
\frac {3\pi}5  , \frac {3\pi}4\ri]\ri\}$. To determine the slopes of the curve $\Phi=0$ we consider $ \arg  ({\d A_+}/{\d \rho})$ modulo $\pi$ (i.e. up to sign).
 Again, we use the fact that an integral with positive measure preserves cones:
 observe that, $\forall s \in [0,1], \rho \in S_L$ the following inequalities hold  
\bea
\arg(\rho+3)<\arg(\rho+3s)<\arg(\rho) <\arg(\rho-3)\ \ \ \Rightarrow \ \ \
- \arg(\rho)<-\arg(\rho+3s)<-\arg(\rho+3)\\
\Rightarrow -\frac {\arg(\rho)}2<\arg \le( \int_0^1 \frac 1 {\sqrt{\rho+3s} } \frac {(1-s) \d s}{ (s+1)\sqrt{s(1-s^2)}}\ri)<-\frac{\arg(\rho+3)}2\;.
\eea
Taking into account \eqref{dAplus} and the above estimates, we find the lower bound
\bea
 \arg\frac {\d A_+}{\d \rho} > -\frac \pi 2 + \overbrace{\frac 1 2 \arg \rho}^{\geq 0} + \overbrace{ {\arg(\rho-3)}}^{\geq \frac {3\pi}5}  
\geq \frac {\pi }{10}\;.
\eea
Likewise, we can obtain an upper bound as follows:
\bea
 \arg\frac {\d A_+}{\d \rho}
< -\frac \pi 2 +\overbrace{\arg(\rho)}^{<\arg(\rho-3)}  \overbrace{-\frac {\arg(\rho+3)}2}^{<0} +\arg(\rho-3) < 
 - \frac \pi 2 +2\overbrace{\arg(\rho-3)}^{\leq \frac {3\pi}4}  \leq  \pi\;.
\eea
The cone $\frac {\pi}{10} <  \arg  ({\d A_+}/{\d \rho}) <\pi$ does not contain the real axis, and hence the tangents of our level curves cannot be horizontal.
Thus
we can parametrize the level curves of $\Phi=0$ within the sector $\mathcal S_L$ by the imaginary part, and point {\bf (5.)} is proved.
 \QED
 
 Assume now that $A\leq B$, the opposite case can be treated symmetrically.  
 Using the  $\Z_5$ action (Lemma \ref{Z5}) we arrange so that the central
 zero $x_2$ lies in the sector $0< {\rm arg} \,x_2 \leq \f{2\pi}{5}$. 
The following lemma describes the range of argument of $x_1$ and $x_3$ under this assumption:
\\
 \noindent
  \begin{minipage}{0.599\textwidth}
   \begin{lemma}
 Let the argument of central zero $x_2$ lie in the range 
  \be
  0 \leq {\rm arg} \, x_2 < \f{2\pi}{5}
  \la{rangex2}\ee
  and the zeros $x_1$ and $x_3$ be labelled such that $A<B$ (i.e. the distance from $x_2$ to $x_1$ is
  greater than the distance from $x_2$ to $x_3$ in the metric $|Q|$). Then the zeros $x_{1,3}$ lie in the following sectors $\Scal_1$ and $\Scal_3$ (see Fig \ref{sectors}):
\be
\f{6\pi}{5}\leq {\rm arg} \, x_1 \leq \f{8\pi}{5}\;,\hskip0.7cm
\f{2\pi}{5}< {\rm arg} \, x_3 \leq \f{4\pi}{5}
\ee
\end{lemma}
\end{minipage} \begin{minipage}{0.4\textwidth}
\begin{center}
\includegraphics{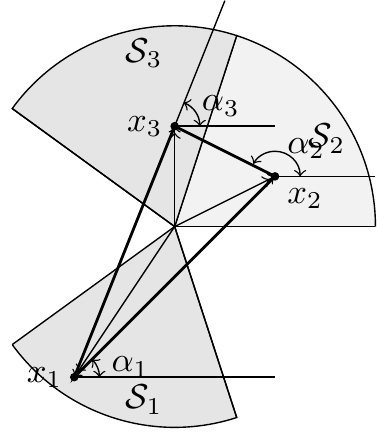}
\end{center}
\captionof{figure}{}
\label{sectors}
\end{minipage}

{\it Proof.} 
The proof uses a continuity argument and the statement of Prop. \ref{ThmW5}.
Using the scaling $\R_+$ action, we shall assume $A+B=1$ throughout this lemma, without loss of generality. 

Consider $A=0, B=1$; then Prop. \ref{ThmW5} implies that $x_3=x_2$ lie on one of the rays $ {\rm e}^{\frac {2i\pi k}5}\R_+$. Using the $\Z_5$ action (Lemma \ref{Z5}) we can assume $k=1$ (see Fig. \ref{rangex2}).  As $A$ increases, the two roots $x_2=x_3$ separate along direction that was computed in Lemma \ref{Lemmaemer} and implies that they move in the two adjacent sectors. We consider the case when the central zero  $x_2$ moves into the sector $(0,\frac {2\pi}5)$; this corresponds to the case "$-$" in \eqref{locscal} (the other case corresponds to a different cell). 

 For small $A$  we thus know that  $x_{1,2,3}$  are in the respective sectors as Fig. \ref{rangex2}; however, as $A$ increases ($A+B=1$), we know from Prop. \ref{ThmW5} that none of the roots can cross the rays $ {\rm e}^{\frac {2i\pi k}5 }\R_+$ unless it is the central zero (which happens only for  $A=B$) , or the curve degenerates (when $A=0$ or $B=0$). Thus, necessarily the three roots $x_{1,2,3}$ remain in the indicated sectors. \QED

\subsection{Variation  of $\arg\Delta_1 $ on $\Qcal_0^\R(-7)$}

The goal of this section is to study the variation of the argument of the modular discriminant $\Delta_1 $ (\ref{moddis}) on the moduli space $\Qcal_0^\R(-7)$.

 The argument of  $\Delta_1 $ (\ref{moddis}) does not change under simultaneous multiplication of all $x_i$ with a positive real constant. Therefore,  $\arg\Delta_1 $ is  constant along any ray passing through the origin in
 the first octant of the $(A,B)$-plane. In particular, $\arg\Delta_1 $ is constant on each 
 of the rays $A=0$, $B=0$ and $A=B$. 
\begin{lemma}\la{lemmavar5}
Let as before the argument of the central zero ${x_2}$ lies between $0$ and $2\pi/5$.
\begin{enumerate}
\item
As $A\to 0$ the arguments of $x_i-x_j$ behave as follows:
\be
{\rm arg} (x_3-x_2)\to\f{13\pi}{20}\;,\hskip0.5cm
{\rm arg} (x_3-x_1)\to\f{2\pi}{5}\;,\hskip0.5cm
{\rm arg} (x_2-x_1)\to\f{2\pi}{5}\;;\hskip0.7cm
{\rm arg} \,\Delta_1 \to \f{29\pi}{10}
\la{limitsA}\ee
\item
As $A= B$ we have
\be
{\rm arg} (x_3-x_2)= \pi-\alpha \;,\hskip0.5cm
{\rm arg} (x_3-x_1)= \f{\pi}{2} \;,\hskip0.5cm
{\rm arg} (x_2-x_1)= \alpha\;;\hskip0.7cm
{\rm arg} \, \Delta_1= 3\pi
\la{limitsAB}\ee
where $\alpha$ is the angle (lying between $0$ and $\pi/2$) formed by the line connecting $x_2$ and $x_3$ with horizontal line. 
\item
As $B\to 0$ the arguments of $x_i-x_j$ behave as follows:
\be
{\rm arg} (x_3-x_2)\to\f{3\pi}{5}\;,\hskip0.5cm
{\rm arg} (x_3-x_1)\to\f{3\pi}{5}\;,\hskip0.5cm
{\rm arg} (x_2-x_1)\to\f{7\pi}{20}\;;\hskip0.7cm
{\rm arg} \,\Delta_1 \to\f{31\pi}{10}
\la{limitsB}\ee

\end{enumerate}

\end{lemma}

{\it Proof.} The limits (\ref{limitsAB}) are trivial: they follow from the symmetry $x_1=\bar{x}_3$ when
$A=B$ and $x_2\in \R$; this proves point $[2]$. 

{\bf [1], [3].} Consider the limits (\ref{limitsA}) arising as $A\to 0$, when $x_2\to x_3$.
When $x_2=x_3$ we have ${\rm arg} \,x_{2,3}=2\pi/5$ and ${\rm arg} \,x_{1}=7\pi/5$; since
all $x_i$ lie on the line going through the origin, we have in this case ${\rm arg} \,(x_{2,3}-x_1)=2\pi/5$.

The angles formed by the tangent directions of $x_{2,3}$ in the process of resolution of a double zero have been computed in Lemma \ref{Lemmaemer}: it follows from it that the limit of ${\rm arg}(x_3-x_2)$ as $x_3\to x_2$ is 

\bea
\lim_{A\to 0} \arg(x_3-x_2) & = &  \frac {2\pi}5+ \frac \pi 4 = \frac {13\pi}{20}
\\
\lim_{A\to 0} \arg(x_3-x_1) & =&  \lim_{A\to 0} \arg(x_2-x_1) = 
 \frac {2\pi}5.
\eea 
and hence
\be
\arg \Delta_1 = \lim_{A\to 0} 2\arg \big((x_2-x_1)(x_3-x_1)(x_3-x_2)\big) = \frac {13\pi}{10} + \frac{8\pi}5=\frac {29\pi}{10}.
\ee 
The proof of \eqref{limitsB} proceeds then in a similar fashion.
\QED

We now are in a position to state the first main analytical result:
\begin{theorem}
\label{propdelta5}
The variation of $\phi_1={\rm arg}\, \Delta_1 $ on the space ${\Qcal}^\R_0(-7)$ from the
 boundary $B=0$ to the boundary $A=0$ equals $\pi/5$. In other words, the monodromy of $\phi_1$ around the point $Q=x^3 dx$
 of ${\Qcal}^\R_0(-7)$, equals $\pi/5$.
 \end{theorem}
{\it Proof.} almost immediately follows from Lemma \ref{lemmavar5}.
To confirm that no additional multiple of $2\pi$ arises when one computes such variation it is sufficient to 
observe that, according to positions of $x_1$ and $x_3$, when ${\rm arg}\, x_2$ is bounded between $0$ and $2\pi/5$ (Fig. (\ref{sectors})) the following inequalities hold:
$$
\pi/5 \leq {\rm arg}(x_3-x_1)\leq 4\pi/5
\qquad 
2\pi/5 \leq {\rm arg}(x_3-x_2)\leq \pi
\qquad 
0 \leq {\rm arg}(x_2-x_1)\leq 3\pi/5
$$
which confirms that the variation of ${\rm arg}\, \Delta_1$ between the boundary $A=0$ and the  diagonal $A=B$ indeed 
equals $\pi/10$; thus the variation of  ${\rm arg}\, \Delta_1$  between the boundaries $A=0$ and $B=0$ equals $\pi/5$. 
\QED

\section{Arguments of  $\Delta^\pm_{-1,-1} $ on the space  $\Qcal_0^\R([-3]^2)$}
\la{space11}

In this section we study  the space $\Qcal_0^\R([-3]^2)$  which represents a local model 
of Kontsevich's boundary $W_{-1,-1} $ of the Strebel combinatorial model of $\Mcal_{g,n}$.  We compute the monodromy of the argument of $\Delta^+_{-1,-1} =x_1^6 x_2^6 (x_1-x_2)^2$ and  $\Delta^-_{-1,-1} =x_1^6 x_2^6 (x_1-x_2)^{26}$  which are the natural
analogs of the modular discriminant on $Q_0^\R([-3]^2)$. The motivation for studying these expressions comes from the theory of Bergman tau-functions $\tau_\pm$ on the spaces of quadratic differentials \cite{contemp, Leipzig}.

\subsection{Space $Q_0([-3]^2)$}
\label{Q33}

An element of the  space $\Qcal_0([-3]^2)$ is a  quadratic differentials $Q$ on the Riemann sphere with  two poles of degree  3 each. The complex  dimension of $\Qcal_0([-3]^2)$ equals 2 and it is  stratified as follows:
\be
\Qcal_0([-3]^2)=\Qcal_0([-3]^2,[1]^2)\sqcup \Qcal_0([-3]^2,2)
\ee

Using a  M\"obius transformation we can assume that the two poles are  $x=0$ and $x=\infty$; the remaining freedom of rescaling of $Q$ 
allows us  to represent it as follows:
\be
Q=\frac{(x-x_1)(x-x_2)}{x^3}(dx)^2
\la{Qgener}
\ee
where $x_1\neq x_2$ are two complex parameters.

Therefore we can identify 
$\Qcal_0([-3]^2,[1]^2)$ with  the space:
\be
\Qcal_0([-3]^2,[1]^2)\simeq \le\{(x_1,x_2\}\in (\C^*)^2,\;\;\; x_1\neq x_2\ri\}/S_2
\ee
For any $Q\in \Qcal_0([-3]^2,[1]^2)$ 
we introduce  the  canonical cover $\Ch$ defined by $v^2=Q$; this is  the elliptic curve
$y^2=x(x-x_1)(x-x_2)$ endowed with the  meromorphic differential 
\be
v=\sqrt{\frac{(x-x_1)(x-x_2)}{x^3}}dx
\ee
having zeros of second order on $\Ch$ at $x_1,x_2$ and poles of second order at $x=0$ and $x=\infty$.

Choosing two canonical cycles $(a,b)$ on $\Ch$ and integrating the differential $v$ over  cycles $a$ and $b$ we get   periods $A$ and $B$ which are defined up to $SL(2,\Z)$ transformation as usual.

According to the next lemma, the periods $(A,B)$ can  be used as local coordinates on  $\Qcal_0([-3]^2,[1]^2)$:
\begin{lemma}
The Jacobian of change of variables from $(x_1,x_2)$ to $(A,B)$ is given by the following formula:
\be
\f{\p(A,B)}{\p (x_1,x_2)}= \pm 2\pi i \frac{x_1-x_2}{x_1 x_2}
\la{Jac11}
\ee
\end{lemma}
{\it Proof.}
The proof is similar to the proof of (\ref{Jacdet}). Namely, to compute the determinant
\be
\frac{\p(A,B)}{\p(x_1,x_2)}
 = \oint_a  \frac {\pa v}{\pa x_1}  \oint_b \frac {\pa v}{\pa x_2} - \oint_b \frac {\pa v}{\pa x_1}    \oint_a  \frac {\pa v}{\pa x_2}
 \la{rbil1}
\ee
with
\be
\frac{\p v}{\p x_1}=-\f{1}{2}\f{1}{x^{3/2}}\left(\f{x-x_2}{x-x_1}\right)^{1/2}\;,\hskip0.7cm
\frac{\p v}{\p x_2}=-\f{1}{2}\f{1}{x^{3/2}}\left(\f{x-x_1}{x-x_2}\right)^{1/2}
\ee
we replace the differential $\p v/\p x_1$ by the holomorphic differential
\be
w= \frac{\p v}{\p x_1}- \f{b}{a}\frac{\p v}{\p x_2}= \f{x_2-x_1}{2 x_1}\f{\d x}{[x(x-x_1)(x-x_2)]^{1/2}}
\ee
(this transformation does not change the determinant).
Then  (\ref{rbil1}) can be computed by applying Riemann bilinear relations to the holomorphic differential $w$ and the meromorphic differential
$ \frac{\p v}{\p x_2}$ which has a pole of order 2 at $x=0$. The computation of residue of $\left(\int^x w\right) \frac{\p v}{\p x_2}$
at $x=0$ in the local parameter $\sqrt{x}$ gives $(x_1 x_2)^{-1}(x_1-x_2)$ which leads to (\ref{Jac11}).
\QED

\subsection{Real slice   $Q_0^\R([-3]^2)$:  Boutroux curves}

By  $Q^\R_0([-3]^2)$ we denote  the space of "Boutroux curves" i.e. the real slice of $Q_0([-3]^2) $ where all
periods  of $v=\sqrt{Q}$
are real. The stratification of  $\Qcal_0^\R([-3]^2)$ is:
$$
\Qcal_0^\R([-3]^2)=\Qcal_0^\R([-3]^2,[1]^2)\cup \Qcal_0^\R([-3]^2,2)\;,
$$
where the  stratum $\Q^\R_0([1]^2,[-3]^2)$ has real dimension 2 and
the stratum  $\Qcal_0^\R([-3]^2,2)$  has real dimension one. The latter corresponds to  differentials 
$v=x^{-3/2}(x-x_1)dx$ such that their period $\int_{x_1^{(1)}}^{x_1^{(2)}} v= 8 x_1^{1/2} $ is real (here $x_1^{(1)}$ and $x_1^{(2)}$ are points on different sheets of Riemann surface of function $x^{1/2}$
having projection $x_1$ on $x$-plane), i.e. $x_1\in \R_+$.

To define periods of $v$ we consider horizontal trajectories of $v$ which look as shown in Fig.\ref{fig2degmid}: two horizontal geodesics always connect $x_1$ and $x_2$ while two other trajectories connect one of the zeros (labelled $x_1$)  to $x=0$ and another zero (labelled $x_2$) to $x=\infty$.  
We label the geodesics  as follows. Let $e_0$ is the horizontal trajectory connecting $x=0$ to $x_1$. When we arrive to $x_1$ following $e_0$  we denote by $e_1$  the first counterclockwise horizontal trajectory issuing from $x_1$ and by  $e_2$ the remaining one; both $e_1$ and $e_2$ end at $x_2$. Denote the length of $e_1$ 
(in the metric $|Q|$) by $A$ and the length of $e_2$ by $B$ (Fig.\ref{fig2degmid}).

\begin{theorem}
The ordered pair of lengths $(A,B)$ defines a one-to-one map between the space ${\Qcal}^\R_0([1]^2,[-3]^2)$ and $\R_+^2$. 
The boundaries $A=0$ and $B=0$ of $\R_+^2$  are identified  and coincide with the space $\Qcal_0^\R([-3]^2,2)$. The remaining period defines a one-to-one map of the space $ \Qcal_0^\R([-3]^2,2)$ to $\R_+$.

The set $A=B$ corresponds to real curves, when $x_{1,2}\in \R$ (the numerical analysis show that for these curves
the approximate value of $J$-invariant equals $J\sim 7791$).
\end{theorem}
{\it Proof.}
To prove that the period map  is invertible on  $ \wh{\Qcal}_0^\R([-3]^2,[1]^2)$ we need to restore the pair $(x_1,x_2)$ 
uniquely knowing periods $A$ and $B$.
The lengths $A$ and $B$ of the critical graph define  a polyhedral surface 
(i.e. surface with flat metric and conical singularities) by gluing 2 half-planes as shown in Fig.\ref{LocalDdeg}. This  uniquely defines the conformal structure i.e 
$J$-invariant of $\Ch$. Moreover, since the labeling of $x_1$ and $x_2$ is fixed, we uniquely reproduce the ratio
$t=x_1/x_2$.

Let us choose the canonical $a$-cycle on $\Ch$ such that $2A=\int_a v$. Then the definition of the period 
$$2A=\int_\a\left[\f{(x-x_1)(x-x_2)}{x^3}\right]^{1/2}\d x$$
can be after the substitution $\tilde{x}=x/x_2$ written as
$$2A=x_2^{1/2}\int_1^t \left[\f{(\tilde{x}-t)(\tilde{x}-1)}{\tilde{x}^3}\right]^{1/2}\d \tilde{x}$$
which fixes $x_2$ uniquely since $t$ is already known.

This proves the isomorphism between $ {\Qcal}_0^\R([-3]^2,[1]^2)$ and the first quadrant $\R_+^2$ in $(A,B)$-plane.
\QED

\begin{figure}[h]
\begin{center}
\includegraphics{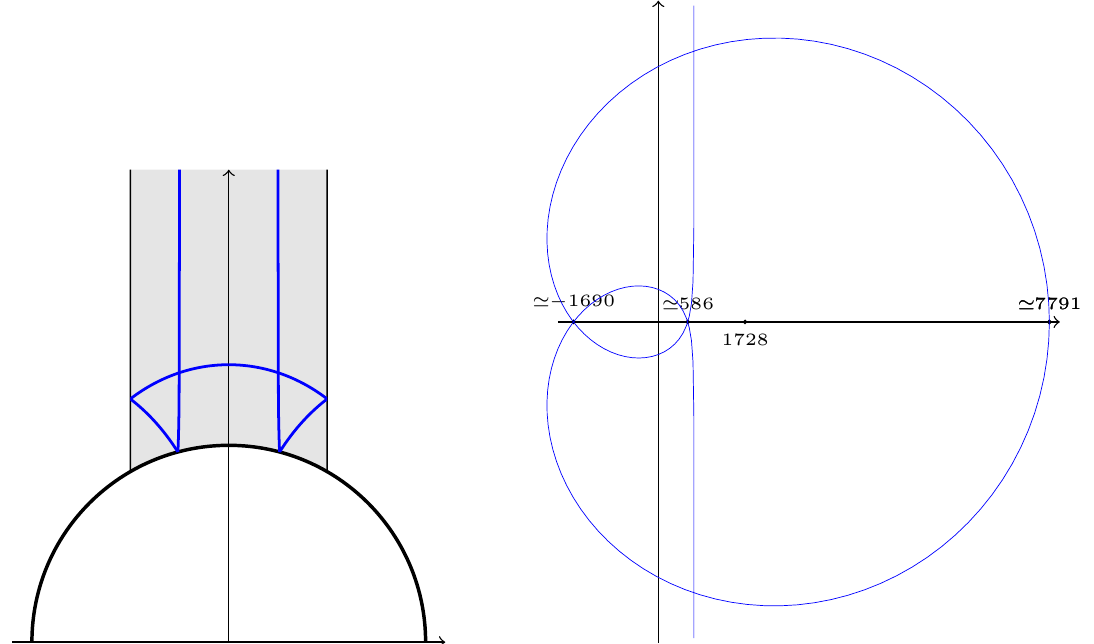}
\caption{Left: the space $\Qcal^\R_0([-3]^2)$ is fibered over the shown curve $\Rcal_{-1,-1}$ in the moduli space of elliptic curves with the fiber $\R_+$.
Right: the set $\Rcal_{-1,-1}$ in the plane of $J$-invariant. Three  real curves correspond to values of $J$-invariant $J_1\simeq -1690$,
$J_2\simeq 586$, $J_3=\simeq 7791$. }
\label{setR11}
\end{center}
\end{figure}

The periods $\sigma$ of canonical coverings $\Ch$ corresponding to points of   $ \Qcal_0^\R([-3]^2,[1]^2)$ 
span  the one-dimensional subset $\Rcal_{1,1}$ of the moduli space $\Mcal_1$ shown in Fig.\ref{setR11}, left. The set 
$\Rcal_{1,1}$ in the plane of $J$-invariant is shown in  Fig.\ref{setR11}, right.
The set $\Rcal_{1,1}$ intersects the real line at the  following theree approximate values of $J$ - invariant:\
\bea
&& J_1\simeq -1690,   \ \ \ (x_1 \simeq 1.5538 + 0.2514i, \ x_2 \simeq 
0.1650 - 0.4893i)
\\
&& J_2\simeq 586.3 \ \ \ (x_1 \simeq 1.2551 + 0.1883i, \ \ x_2\simeq 0.6865 - 0.2986i)  \\
&& J_3\simeq 7791,\ \  \ (x_1\simeq  1.8037, \ \ x_2 = -0.3797)
\eea

All three values $J_i$ correspond to curves $\Ch$ possessing a real involution. However, only the value $J_3$ corresponds to
the Boutroux-Krichever curve where the Abelian differential $v=\sqrt{Q} $ is invariant under this involution i.e.
$v(\bar {x})=\overline{v(x)}$.

   \subsection{Variation of $\arg \Delta^\pm_{-1,-1} $ on $\Qcal^\R_0([-3]^2)$}

The goal of this section is to study the variation of the arguments of the expressions
\be
 \Delta^+ _{-1,-1} := x_1^6 x_2^6 (x_1-x_2)^2,\  \ \qquad  \Delta^- _{-1,-1} := x_1^6 x_2^6 (x_1-x_2)^{26}
\la{defDet}\ee
over the space $ \Qcal_0^\R([-3]^2) $. The origin of  the expressions 
(\ref{defDet}) lies in the formalism of Bergman tau-function \cite{Annalen,JDG,MRL,contemp}.

As well as the previos case, the arguments of  $\Delta^\pm_{-1,-1} $ (\ref{defDet}) do not change under simultaneous multiplication of all $x_i$ with positive real constant. Therefore,  $\arg\Delta^\pm_{-1,-1} $ is  constant along any ray passing through the origin in
 the first octant of the $(A,B)$-plane. In particular, $\arg\Delta^\pm_{-1,-1} $ is constant on each 
 of the rays $A=0$, $B=0$ and $A=B$. 

\begin{proposition}
\label{W331}
Suppose that $Q =\frac {(x-x_1)(x-x_2)}{x^3} \d x^2\in \mathcal Q_0^\R([-3]^2)$. 
Then
\begin{enumerate}
\item  one root $x_i$ is real if and only if the other is: moreover, they either coincide ($x_1=x_2 >0$) or one of them (say, $x_2$)  is positive and the other ($x_1$)  is negative. 
In this latter case, the horizontal trajectories connect $x_1$ to $0$ along the negative axis, and $x_2$ to $+\infty$ along the positive axis. 
\item if one root is in the upper half plane then the other is in the lower half plane;
\end{enumerate}
\end{proposition}
The proof of this  proposition is technical; it is contained in  the Appendix. In Fig. \ref{figdeg2} we show how the branch points $x_i$ and the horizontal trajectories evolve as periods $(A,B)$ move from the point $(1,0)$ to the point $(0,1)$ along the straight line.

The previous  proposition implies

\begin{proposition}
\label{proparg}
\begin{enumerate}
\item
As $A\to 0$ we have
\be
{\rm arg}\, x_{1,2}\to 0   \hskip0.7cm
{\rm arg}\, (x_2-x_1)\to\, \pi/4
\ee

\item
As $A=B$ we have 
\be
{\rm arg}\, x_1=  0   \hskip0.7cm
{\rm arg}\, x_2= \, \pi  \hskip0.7cm
{\rm arg}\, (x_2-x_1) = \,0  
\ee

\item

As $B\to 0$ we have
\be
{\rm arg}\, x_1\to \,2\pi  \hskip0.7cm
{\rm arg}\, x_2\to\, 0  \hskip0.7cm
{\rm arg}\, (x_2-x_1)\to\, - \pi/4
\ee
\end{enumerate}
\end{proposition}
{\it Proof.} {\bf [1]}. As the curve $\Ch$ degenerates, item $1$ of Prop. \ref{W331} implies that $\arg_{1,2} \to 0$ and $x_{1,2} \to x_0>0$; in this limiting case the flat coordinate is $z(x) = \frac {\sqrt{x}}2+ \frac {x_0}{2\sqrt{x}}$ and the critical trajectories $\Re z=0$ are easily seen to consist of $\R_+ \cup \{|x|=x_0\}$. 

By Prop. \ref{W331} the two roots remain in opposite half-planes; with our choice of $A,B$, it is clear that as $A $ decreases, $x_1$ moves in the lower half-plane from $\R_-$  and $x_2$ in the upper half-plane from $\R_+$.

When $A\to 0$, we know that the curve degenerates and both roots converge to a point on $\R_+$; since they are confined in their respective half-planes, the statement about the increments of $\arg x_{1,2}$ follows immediately.  As $A\to 0_+$, an analysis entirely similar to that of Lemma \ref{Lemmaemer}  (which we do not repeat), shows that $\lim_{A\to 0_+}\arg(x_2-x_1) = \frac \pi 4$. See also the numerically accurate Fig. \ref{figdeg2}. 

{\bf [2].} 
From part $1$ of Prop. \ref{W331} it follows that for $A=B$ we have  $x_1\in \R_-, x_2\in \R_+$ and hence $\arg(x_2-x_1)=0$. 

{\bf [3].} This  case is obtained from case [1] by sending $v(x)\mapsto \ov {v(\ov x)}$. 
\QED

Finally we are in the position to formulate the second main analytical result: 
\begin{theorem}
\label{prop2}
The monodromies  of ${\rm arg}\,\Delta^\pm _{-1,-1} $ \eqref{defDet} on the space ${\Qcal}_0^\R([-3]^2)$ 
along  an elementary non-contractible loop equal $13\pi$ and $25\pi$, respectively.
\end{theorem}
{\bf Proof.}
It suffices to refer to Prop. \ref{proparg} to compute the increments of ${\rm arg}\,\Delta^\pm _{-1,-1} $ 
from the diagonal $A=B$ to the boundary $A=0$. Namely, for $A=B$ we have $\arg x_1 = -\pi$, $\arg x_2=0$, $\arg(x_2-x_1)=0$ and for $A=0$ we have (in the limit) $\arg x_{1,2}=0$, $\arg(x_2-x_1)=\frac \pi 4$. Hence the variations are $\f{\pi}{2} + 6\pi = \frac {13}2\pi$ for $\Delta^+_{-1,-1}$ and $
\frac {13}2 \pi +6\pi = \frac {25}2 \pi$ for $\Delta^-_{-1,-1}$.  To get the monodromy over the whole cell we must duplicate these increments.
\QED

\begin{figure}

\resizebox{0.75\textwidth}{!}{
\begin{minipage}{1\textwidth}
\includegraphics{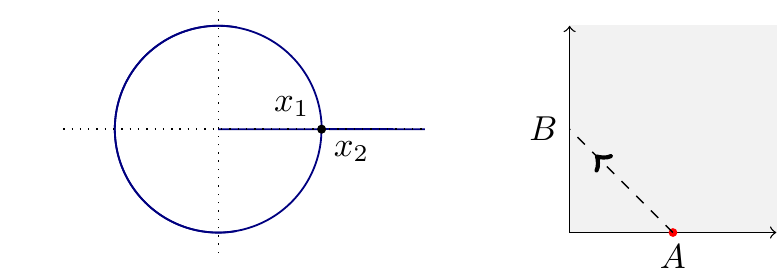}
\includegraphics{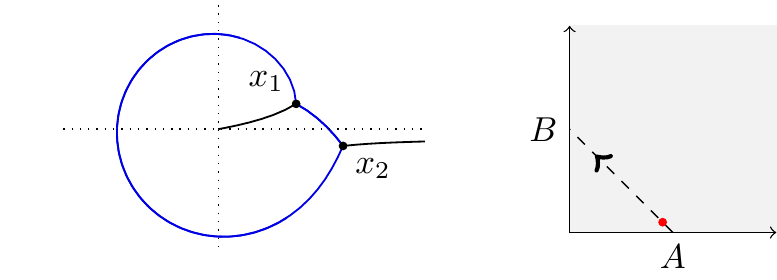}
\includegraphics{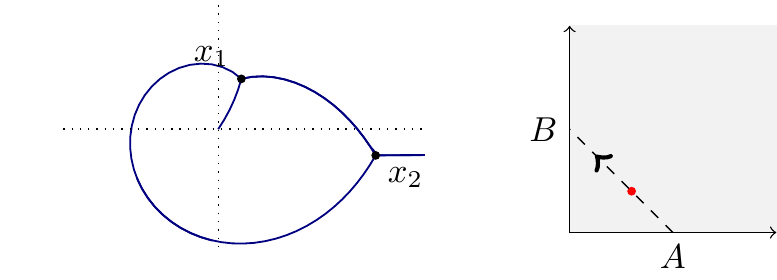}
\includegraphics{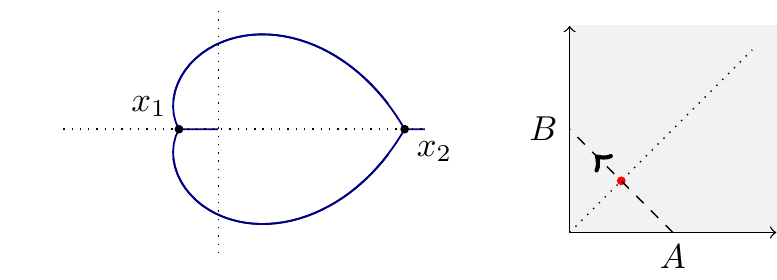}
\includegraphics{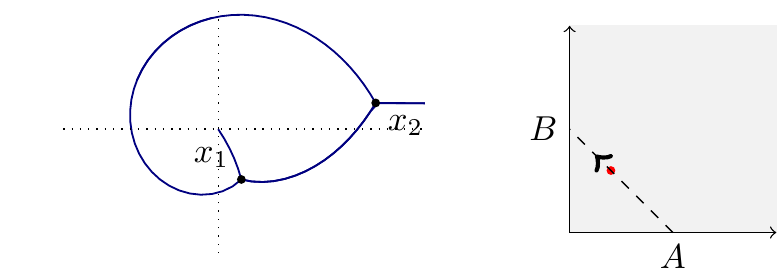}
\includegraphics{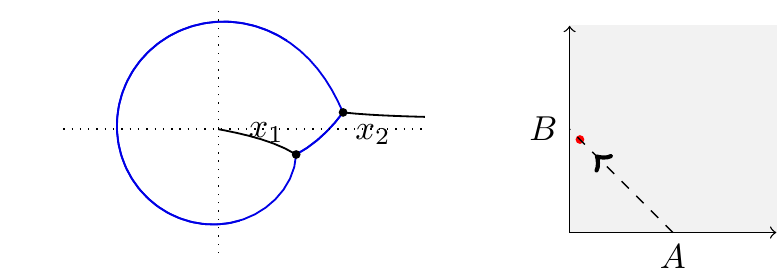}
\includegraphics{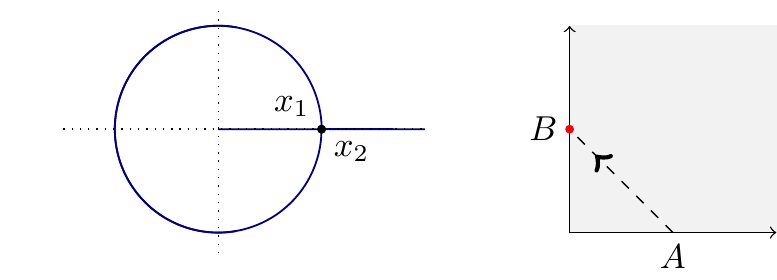}
\end{minipage}
}
\captionof{figure}{ As the periods $(A,B)$ move from the point $(1,0)$ to the point $(0,1)$ along the line $A+B=1$,  the root marked $x_1$ rotates counterclockwise around the origin. This move corresponds to monodromy around the point
$Q=(dx)^2/x$  on the space    $\wh{Q}_0^\R([-3]^2)$.}
\label{figdeg2} 
\end{figure}

\paragraph {\bf Acknowledgements.} D.K. thanks Peter Zograf for interesting discussions. The research of M. B.   was supported in part by the Natural Sciences and Engineering Research Council of Canada grant
RGPIN/261229--2011. The research of D.K.   was supported in part by the Natural Sciences and Engineering Research Council of Canada grant
RGPIN/3827-2015,  Alexander von Humboldt Stiftung and GNFM Gruppo Nazionale di Fisica Matematica. Both authors were supported by  the FQRNT grant "Matrices Al\'eatoires, Processus Stochastiques et Syst\`emes Int\'egrables" (2013--PR--166790).
D.K.   thanks  the International School of Advanced Studies (SISSA) in Trieste and Max-Planck Institute for Gravitational Physics in Golm (Albert Einstein Institute) for their hospitality and support during the preparation of this paper.

\appendix
\section { Proof of Proposition  \ref{W331}} 
{[\bf 1]} The Boutroux property (i.e. the reality of all periods) is invariant under the action of $\R_+$  that maps $x\mapsto \lambda x, \ x_j \to \lambda x_j,\ Q \mapsto \lambda Q$. Suppose that one of the roots $x_j$ is real; depending on its sign, using the above scaling, we can restrict ourselves by the case that this root equals $\pm 1$. 
Let us therefore study  the differentials 
\be
v_+ = \sqrt{\frac { (x-q)(x-1)}z }\frac{\d x}{x},\ \hbox { and } \ v_- = \sqrt{\frac { (x-q)(x+1)}x }\frac{\d x}{x}\;.
\ee
\paragraph{Degenerate case.} 
Suppose $q=1$ for $v_+$ or $q=1$ for $v_-$; then the resulting differential  can be integrated  explicitly 
\be
\int^x_1 v_+ = 2\le(\sqrt{x}+\frac 1{\sqrt{x}} -2\ri)\\;,  \qquad  \int^x_{-1} v_- =2\le(\sqrt{x} - \frac 1{\sqrt{x}}-2i\ri)\;.
\ee
In the first case we see easily that the horizontal trajectories (where the real part of the integral vanishes) are $|x|=1$ and $\R_+$. 
In the second case there is no trajectory connecting $-1$ to itself; this is seen by computing the integral along closed path connecting $x=-1$ to itself around the origin, and noticing that the result is actually imaginary. 

\paragraph{Nondegenerate case.} Let us show that the Boutroux condition makes $q$ to be real and of the opposite sign (i.e. $q$ is negative in the case of $v_+$ and positive in the case of $v_-$).
We consider two cases separately.

\noindent \begin{minipage}{0.74\textwidth}\paragraph{Case $+$.}
Let us rewrite the periods by an integration by parts:\\
\bea
 && H (q):= \frac 1 2\int_{\gamma_1}  \sqrt{\frac {(x-  1)(x-q)}x}\frac {\d x}x \mathop{=}^{b.p.}
\int_0^1    \frac {q+1-2x}{\sqrt{x(1-x)(q-x)} }{\d x}\;,
\label{H}
\\
&&
J(q) := \frac 1 2\int_{\gamma_q}  \sqrt{\frac {(x-1)(x-q)}x}\frac {\d x}x \mathop{=}^{b.p.\atop x \mapsto qx}
 \int^{1}_0    \frac {1+(1 -2x)q}{\sqrt{x( 1-x)(1-qx)} }{\d x}\;.
\label{J}
\eea
\end{minipage}\,\,\ 
\begin{minipage}{0.25\textwidth}
\begin{center}
\resizebox{1\textwidth}{!}
{
\includegraphics{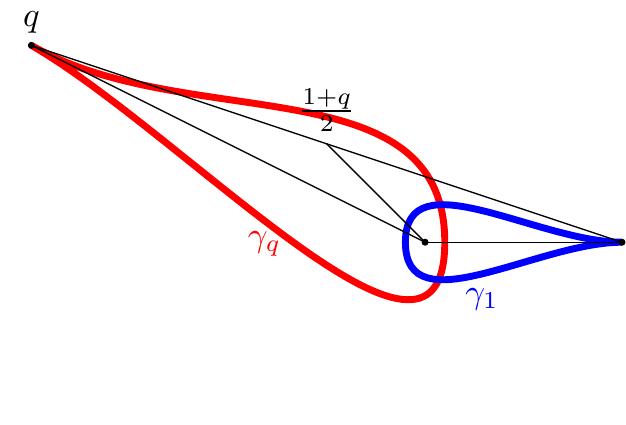}
}
\end{center}
\captionof{figure}{The contours of integration $\gamma_1, \gamma_q$.}
\end{minipage}\\

Our goal is to show that 
the equations $\Im I(q) = 0 =\Im J(q)$ have only one  solution $q^\star\in \R_-$.

We start by showing that the locus of the equation $\Im J=0$ is the set  $q\in (-\infty, 1]$.
The function $J(q)$ is manifestly analytic in the domain $\mathcal D = q\in \C\setminus [1,\infty)$ and has a  branch point at $q=1$.
A simple analysis of the integral \eqref{J}  shows that 
\begin {enumerate}
\item $\ds \lim_{q\to 0 + i0} J(q) =  \int_0^1 \frac {\d x}{\sqrt{x(1-x)} }= \pi$. 
\item $ \ds \lim_{q\to 1+ i0} J(q) = 4$;
\item $J( (-\infty, 1] + i0) = (-\infty, 4]$;
\item for $q>1$ the imaginary part of $J(q)$ (boundary value from upper half plane) is positive;  this is the estimate  in \eqref{1122} below.
\end{enumerate}

 On the domain $\mathcal D$ we also have the Schwartz symmetry  $\ov {J(\ov q)}=J(q)$ and hence it suffices to study it in the upper half plane. 

 The locus $\Im J=0 $ is  the union of possibly several smooth arcs, each of which must be extensible indefinitely until they reach either a singularity of $J(q)$ or a zero of $J'(q)$. We now show that $J'(q)\neq 0$ in $\mathcal D$; this is easily seen by explicit computations as follows
\be
J'(q) = -\frac 1 4 \int_{\gamma_q} \frac {\d x}x \sqrt{ \frac { x-1}{x(x-q)}}  \mathop{=}^{b.p.} 
-\frac 1 {4q}  \int_{\gamma_q}  \sqrt{\frac {q-x}{ x(1-x)}} \d x
\mathop{=}^{z=qt} -\frac {1 } 2 \int_0^1 \sqrt{\frac {1-t}{t(1-qt)}} \d t
\ee
The last integral is clearly always nonzero for all $q\in\mathcal D$.

Since there are no zeroes of $J'(q)$ in $\mathcal D$, any branch of $J^{-1}(\R)$ that may be in $\mathbb H$ must extend to infinity along some direction or be a bounded curve starting and ending at $q=1$; this latter case is excluded because $\Im J$ is harmonic and  bounded at $q=1$ and there are no singularities in $\mathcal D$. 

Note that as $|q|\to \infty$  in the upper half plane, we have from \eqref{J} that  $J(q) \sim i\sqrt q  \ln (1-q) \sim \sqrt{-q} \ln (-q)$ and the real part (to leading order in $|q|$) is easily seen not to change sign for $0<\arg(q)<\pi$. Hence there are no branches of $J^{-1}(\R)$ extending to infinity. 
%
%
To prove item $4$ from the above list we observe that 
for $q>1$
\bea
\Im J(q) = \sqrt{q} \int_{\frac 1 q} ^1  \frac {q^{-1} +1-2x}{\sqrt{x(1-x)(x-q^{-1})}}
=\\
=
\sqrt{q}  \int_{q^{-1}}^{\frac {q^{-1}+1}2}  \frac {\overbrace{q^{-1} +1-2x}^{>0}}{\sqrt{x(1-x)(x-q^{-1})}} \d z
+
\sqrt{q}  \int^1_{\frac {q^{-1}+1}2}  \frac {\overbrace{q^{-1} +1-2x}^{<0}}{\sqrt{x(1-x)(x-q^{-1})}}\d x
\label{1122}
\eea
To see that the result is positive we proceed to a simple estimate of the two integrals:
\bea
\label{J1}
\int_{q^{-1}}^{\frac {q^{-1}+1}2}  \frac {{q^{-1} +1-2z}}{\sqrt{x(1-x)(x-q^{-1})}} \d x 
>
 \int_{q^{-1}}^{\frac {q^{-1}+1}2}  \frac {{q^{-1} +1-2x}}{\sqrt{\frac {q^{-1}+1}2(1-x)(x-q^{-1})}} \d x
\eea
and
\bea
\label{J2}
 \int^1_{\frac {q^{-1}+1}2}  \frac {q^{-1} +1-2x}{\sqrt{x(1-x)(x-q^{-1})}}\d z
>
  \int^1_{\frac {q^{-1}+1}2}  \frac {q^{-1} +1-2x}{\sqrt{\frac {q^{-1}+1}2(1-x)(x-q^{-1})}}\d x
  \mathop{=}^{x\mapsto 1+q^{-1}-x}\\
  =
  \int^{q^{-1}}_{\frac {q^{-1}+1}2}  \frac {q^{-1} +1-2x}{\sqrt{\frac {q^{-1}+1}2(1-x)(x-q^{-1})}}\d x
 =  -\int_{q^{-1}}^{\frac {q+1}2}  \frac {q^{-1}+1-2x}{\sqrt{\frac {q^{-1}+1}2(1-x)(x-q^{-1})}}\d x
\eea
and therefore the sum \eqref{1122} is strictly positive.
In summary, the only branch $\Im J=0$ is $(-\infty, 1]$. 
This ends the proof that $J(q)$ is real only on  $(-\infty, 1]$.
\vskip 0.5cm
We now analyze the other integral $H(q)$ \eqref{H} on the real $q$--axis and show that it
%
takes real values  only for   $q\in [1,\infty) \cup \{q^\star\}$ where  $q^\star <0$. 
For  $q\geq 1$ the integrand is  a real valued function for $x\in (0,1)$.

Let us show that for  $q\in [0,1)$  the imaginary part of $H(q)$  is strictly positive; indeed, for $0<x<q$  the integrand is real--valued and so the imaginary part is (this is analogous to the estimates \ref{J1}, \ref{J2} with the replacement $q^{-1} \mapsto q$)
\be
\Im H(q) = \int_q^1  \frac {q +1-2x}{\sqrt{x(1-x)(x-q)}}
= \int_q^{\frac {q+1}2}  \frac {\overbrace{q +1-2x}^{>0}}{\sqrt{x(1-x)(x-q)}} \d z
+
 \int^1_{\frac {q+1}2}  \frac {\overbrace{q +1-2x}^{<0}}{\sqrt{x(1-x)(x-q)}}\d z
\label{112}
\ee
To see that the result is positive we proceed to a simple estimate of the two integrals:
\bea
\int_q^{\frac {q+1}2}  \frac {{q +1-2x}}{\sqrt{x(1-x)(x-q)}} \d x 
>
 \int_q^{\frac {q+1}2}  \frac {{q +1-2x}}{\sqrt{\frac {q+1}2(1-x)(x-q)}} \d x
\\
 \int^1_{\frac {q+1}2}  \frac {q +1-2x}{\sqrt{x(1-x)(x-q)}}\d x
>
  \int^1_{\frac {q+1}2}  \frac {q +1-2x}{\sqrt{\frac {q+1}2(1-x)(x-q)}}\d x
  \mathop{=}^{z\mapsto 1+q-x}
    -\int_q^{\frac {q+1}2}  \frac {q +1-2x}{\sqrt{\frac {q+1}2(1-x)(x-q)}}\d x
\eea
and therefore the sum \eqref{112} is strictly positive and  so $\Im H(q)>0$ for $q\in [0,1)$.

Now let $q<0$; then 
\be
H(q) =  i \int_0^1 \frac {x + |q| + (x -1)}{\sqrt{x(1-x)(|q| + x)}}\d x
\ee
so that it is purely imaginary. The imaginary part is a monotone function:
\be
\frac {\d}{\d q} \Im H(q) = - \frac 1 2\int_0^1 \frac {|q| + 1}{\sqrt{x(1-x)(|q| + x)} ( x+|q|)}\d x <0.
\ee
It is easily seen that  $\lim_{q\to -\infty}H(q) = i\int_0^1 \frac {\d x}{\sqrt{x(1-x)}} = i\pi$ while as $q\to 0_-$ the value of $H(q)$ diverges logarithmically to $-\infty$. Therefore it vanishes at a single point $q^\star<0$.

Numerical analysis shows that
 $$q^\star \simeq -0.21048557\;.$$

\paragraph{Case $-$.}
Two independent homological coordinates are given by the integrals
\bea
&& 
\wh H (q):= \frac 1 2\int_{\gamma_{-1}}  \sqrt{\frac {(x+  1)(x-q)}x}\frac {\d x}x \mathop{=}^{b.p.\atop x\mapsto -x}
\int_0^{1}    \frac {1-q+2x }{\sqrt{x(1-x)(q+x)} }{\d x}\;,
\\
&& \wh J(q) := \frac 1 2\int_{\gamma_q}  \sqrt{\frac {(x+1)(x-q)}x}\frac {\d x}x \mathop{=}^{b.p.\atop x\mapsto qx}
i \sqrt {q}\int^{1}_0    \frac {1-q^{-1} -2x}{\sqrt{x( x + q^{-1})(1-x)} }{\d x}\;.
\eea
\begin{itemize}
\item We first show that  $\Im q>0$ implies   $\Im \wh H(q)<0$ and that $\wh H(q)\in \R$  only for $q>0$; this is easily seen by rewriting  $\wh H(q) $ as 
\bea
\wh H(q) 
  =  \int_0^{1}    \frac {1+3x }{\sqrt{x(1-x)(q+x)} }{\d x} +   \int_0^{1}    \frac {-\sqrt{q+x }}{\sqrt{x(1-x)} }{\d x} \;.
\label{whH}
\eea
Suppose $\Im q>0$;  a simple inspection of phases confirms that the integrands in  both  integrals have negative imaginary part. 

Therefore the requirement $\Im \wh H=0$ implies $q\in \R$.
Suppose that  $q<0$;
then  the integrals in \eqref{whH}  (the boundary value being taken from $\Im q>0$) acquire an imaginary contribution from the integration over $x\in [0, \min(1,-q)]$ and both contributions have negative imaginary part.  We conclude that $q\in \R_+$. 

\item 
Having now restricted $q$ to $\R_+$ we must see when the integral $\wh J(q)$ is real for $q\in \R_+$.
For $q>0$ it is clear that $\wh J(q) \in i\R$ and hence we are looking at the zeroes of $\wh J(q)$ on $q\in \R_+$. 
An elementary analysis shows  that $\Im \wh J(q)$  is an increasing function of $q\in \R_+$ which changes sign only once. Hence $\wh J(q)$ is real only at one point $\tilde{q}r$, which must then be obtained by  the dilation action from the previous case, so that
 $\tilde{q} = -\frac 1{q^\star}>0$. 
\end{itemize}
We finally show that the negative root is connected to zero by a horizontal trajectory; since the whole set of trajectories trajectories must be invariant under $z \mapsto \ov z$, if a trajectory issuing from either root runs along the real axis, it must be either a segment to $0$ or a ray to $\infty$. An elementary local analysis shows that the asymptotic directions of every critical trajectory falling towards $0$ go along $\R_-$ and those approaching  $\infty$ go along $\arg(z)=0$. Thus necessarily the negative root is connected to $0$ and the positive root is connected to $+\infty$.

{\bf [2]} In item {\bf [1]} we proved  that either both roots have nonzero imaginary part or they are both real; 
therefore, if they are not real, they either are in the same upper/lower half plane or in the opposite ones. 
Let us show it is the latter case. To see that it is sufficient to compute the differential with respect to the periods at the point where both roots are real. 

So let $v = \sqrt{ (x-x_1) (x-x_2)/x}\,\d x$  with $x_1=1, x_2 = q^\star$; we want to find what happens under an infinitesimal Boutroux deformation of $x_0,x_1$ (a deformation that preserves the Boutroux condition). A trivial computation yields
(dot denotes a derivative with respect to some real deformation parameter $t$ at $t=0$)
\be
\dot v  = \frac {-x_1 \dot x_2 -x_2 \dot x_1+ ( \dot x_1 + \dot x_2) x }{2x\sqrt{x(x-x_1)(x-x_2)}}\d x
\ee
where $\dot x_1, \dot x_2$ are constrained by the requirement that all periods are real: for example if the deformation is $\dot A =1, \dot B=0$ then we would require that the periods of $\dot v$ are the given ones, which uniquely determines the coefficients $\dot x_1,\dot x_2$. Since we are computing the deformation at $x_1=1, x_2=q^\star\in \R_-$ we can write, 
\be
\dot v  = \frac { -\dot x_2 - q^\star \dot x_1 + ( \dot x_1+ \dot x_2) x }{2x\sqrt{x(x-1)(x-q^\star)}}\d x
\ee
Consider now the period from $x=1$ to $x=\infty$;
\bea
-(\dot x_2 + q^\star \dot x_1) C_0+ (\dot x_1+\dot x_2) C_1 =\alpha\in \R\;;
\\
C_j = \int_{1}^\infty \frac {x^j\d x}{\sqrt{x(x-1)(x-q^\star)}} .\ \ \ \label{sys1}
\eea
It is evident that $C_j$'s are real and of the same sign depending on the chosen determination of the square root (say, positive). Moreover it is also evident that $|C_1|>|C_0|$ because the integration is on $[1,\infty)$; taking the imaginary part of  \eqref{sys1} 
yields 
\be
-C_1 \Im  ( q^\star \dot x_1   + \dot x_2)+ C_0 \Im (\dot x_1+\dot x_2) = 0\in \R
\ \ \Rightarrow  \ \  
\Im \dot x_2 =\frac {C_0- q^\star C_1}{C_0-C_1} \Im \dot x_1.
\ee
The coefficient $\frac {C_0- q^\star C_1}{C_0-C_1}  $  is negative  because  $C_0,C_1$  are of the same sign,  $q^\star \simeq -0.2$ is negative and $|C_1|>|C_0|$. Thus the initial directions of the motion of the two roots point to the opposite half-planes.  By a simple continuity argument, the points $x_{1,2}$  cannot cross  the real axis without falling onto the case covered by the previous analysis. Hence they remain confined within opposite half-planes.
\QED

\end{document}